\documentclass[preprint,11pt]{elsarticle}

\usepackage{amsmath,amssymb,amsfonts}
\usepackage{graphicx,graphics}
\usepackage{subfigure}
\usepackage{multirow}
\usepackage{pgfplots}
\usetikzlibrary{plotmarks}
\usepackage{hyperref}

\usepackage{caption}

\usepackage{longtable}
\usepackage{rotating}

\usepackage{geometry}
\usepackage{listings}
\usepackage{color}

\definecolor{dkgreen}{rgb}{0,0.6,0}
\definecolor{gray}{rgb}{0.5,0.5,0.5}
\definecolor{mauve}{rgb}{0.58,0,0.82}

\usepackage{amsmath,amsfonts,amssymb}
\usepackage{subfigure}
\usepackage{multirow}
\usepackage{lscape}
\usepackage{setspace}

\newcommand{\aaa}{\mathbf{a}}

\newcommand{\bb}{\mathbf{b}}

\newcommand{\dd}{\mathbf{D}}

\newcommand{\qq}{\mathbf{q}}		
\newcommand{\uu}{\mathbf{u}}
\newcommand{\xx}{\mathbf{x}}
\newcommand{\cn}{\mathbf{n}}
\newcommand{\vv}{\mathbf{v}}

\newcommand{\bvsig}{\boldsymbol{\sigma}}
\newcommand{\bveps}{\boldsymbol{\varepsilon}}

\newcommand{\Eref}[1]{Equation (\ref{#1})}
\newcommand{\fref}[1]{Figure \ref{#1}}
\newcommand{\Erefs}[1]{Equations (\ref{#1})}

\newtheorem{remark}{Remark}

\begin{document}


\begin{frontmatter}



\title{Numerical evaluation of stress intensity factors and T-stress for interfacial cracks and cracks terminating at the interface without asymptotic enrichment}

\author[a]{Sundararajan Natarajan\fnref{label1}\corref{cor1}}
\author[a]{Chongmin Song}
\author[b]{Salim Belouettar}

\address [a]{School of Civil and Environmental Engineering, The University of New South Wales, Sydney, NSW 2052, Australia.}
\address[b]{Centre de Recherche Public Henri Tudor, 29, Avenue John F Kennedy, L-1855 Luxembourg-Kirchberg, Luxembourg.}

 \fntext[1]{School of Civil and Environmental Engineering, The University of New South Wales, Sydney, NSW 2052, Australia. Tel: +61(2)9385 5030. Email: s.natarajan@unsw.edu.au; sundararajan.natarajan@gmail.com}

\begin{abstract}
In this paper, we extend the recently proposed extended scaled boundary finite element method (xSBFEM)~\cite{natarajansong2013} to study fracture parameters of interfacial cracks and cracks terminating at the interface. The approach is also applied to crack growth along the interface and crack deflecting into the material within the context of linear elastic fracture mechanics. Apart from the stress intensity factors, the T-stress can be computed directly from the definitions, without any requirement of path independent integrals. The method aims at improving the capability of the extended finite element method in treating crack tip singularities of cracks at interfaces. An optimum size of the scaled boundary region is presented for multimaterial junctions. The proposed method: (1) does not require special numerical integration technique; (2) does not require a priori knowledge of the asymptotic fields and (3) the stiffness of the region containing the crack tip is computed directly. The robustness of the proposed approach is demonstrated with a few examples in the context of linear elastic fracture mechanics. A discussion on the crack growth along the interface and crack deflecting into the material is also presented.
\end{abstract}


\end{frontmatter}

\section{Introduction}

The finite element method (FEM) reformulates the continuous boundary and initial value problems into equivalent variational forms. The FEM requires the domain to be subdivided into non-overlapping regions called `elements' and local polynomial representation is used for the fields within the element. The FEM with piecewise polynomials are inefficient to deal with internal discontinuities such as material interfaces or singularities. The material interfaces can be captured by a conforming mesh and the singularities can be captured by using quarter-point elements~\cite{barsoum1977} or by locally enriching the finite element approximation space with singular functions~\cite{strangfix1973}. In case of the quarter point elements, the singularity is captured by moving the element's mid-side node to the position one quarter of the way from the crack tip. This was considered to be a major milestone in applying the FEM for linear elastic fracture mechanics (LEFM). Tracey~\cite{tracey1971} and Atluri \textit{et al.,}~\cite{atlurikobayashi1975} proposed a new element that has an inverse square root singularity near the crack. The main advantage of this element is that the stress intensity factors (SIFs) can be computed more accurately. Strang and Fix in their well known book~\cite{strangfix1973}, proposed to add singular functions to approximate the displacement field near the singularity. Benzley~\cite{benzley1974} developed an arbitrary quadrilateral element with singular corner node by `enriching' a bilinear quadrilateral element with singular terms. All the above mentioned approaches still relies on a conforming finite element mesh. In an effort to overcome the limitations of the FEM, meshfree methods were introduced. Treatment of evolving discontinuities in meshfree methods is more straightforward because it does not require conforming mesh or mesh adaptation as the discontinuities evolve~\cite{rabczukgracie2010,rabczukbelytschko2006,rabczukzi2008a,zirabczuk2007,bacnguyen-xuan2013,chau-dinhzi2012}. Three dimensional crack propagation have been presented in~\cite{rabczukbelytschko2007,zhuangaugarde2012,rabczukzi2010,bordasgzi2008,rabczukzi2008a,rabczukbordas2007}.On another albeit related front, in the past decade reducing the meshing burden has spurred interest in the research community, whilst still retaining the finite element framework~\cite{melenkbabuvska1996,belytschkoblack1999,garciafancello2000,rabczukbelytschko2004,rabczukzi2008,nguyen-xuanliu2013}. Of particular note, in this study, we focus our attention on enrichment methods, especially, on the methods which are built on the finite element framework, for example, the extended FEM (XFEM). However, the idea presented here can be coupled with other methods, such as the meshless methods~\cite{belytschkoorgan1995} and the boundary element method~\cite{bird2012} to name a few.

Based on the seminal work of Babu\v{s}ka~\textit{et al.,}~\cite{melenkbabuvska1996}, Belytschko's group in 1999~\cite{belytschkoblack1999} introduced the XFEM that can model crack propagation and strong discontinuities with minimal remeshing. Since its inception, the XFEM has emerged as a versatile tool to handle internal discontinuous surfaces or moving boundaries independent of the underlying finite element discretization. Built on a finite element framework, the method relies on a modification (augmentation, enrichment) of the finite element spaces to capture the internal discontinuities. The success of the XFEM when applied to moving boundary problems, esp, crack growth, relies on the a priori knowledge of the set of functions that span the asymptotic fields ahead of the crack tip. The requirement for a priori knowledge of the asymptotic fields hinders the application of the XFEM directly to heterogenous materials for which the asymptotic fields do not exist in closed form or are very complex. When the required functions are not known a priori, these functions can be computed numerically~\cite{duartekim2008,waismanbelytschko2008,menkbordas2010,mousavigrinspun2011,mousavigrinspun2011a}. Nevertheless, all the approaches adopted in the literature depend on the following two steps:
\begin{itemize}
\item If the enrichment functions are unknown a priori, the enrichment functions are computed numerically, for example, the work of Menk and Bordas~\cite{menkbordas2010}, Mousavi \textit{et al.,}~\cite{mousavigrinspun2011} and Zhu~\cite{zhu2012}.
\item The numerically computed functions are then used as enrichment functions within the XFEM/GFEM framework.
\end{itemize}
In addition to the requirement of the knowledge of asymptotic fields, augmenting these functions to the FE approximation basis leads to other difficulties such as, blending problem, numerical integration of enrichment functions, ill-conditioning and additional unknowns. Recent research has been directed towards improving the method by alleviating some of the aforementioned difficulties (for example, hybrid crack element~\cite{xiaokarihaloo2007}, spider XFEM~\cite{chahinelaborde2008}, hybrid analytical XFEM~\cite{rethoreroux2010}, stable GFEM~\cite{babuvskabanerjee2012}, etc.,) whilst the other focus has been in developing new enrichment functions when XFEM is applied to heterogeneous materials (for example, enrichment functions for orthotropic materials~\cite{asadpouremohammadi2007,asharimohammadi2011,hattorirojas-diaz2012}, for functionally graded materials~\cite{dolbowgosz2002}, for crack at bimaterial interfaces~\cite{sukumarhuang2004,liuxiao2004} and cracks terminating at the interface~\cite{bouhalashao2013}, to name a few). However, in spite of these advancements, to accurately evaluate the fracture parameters (such as the stress intensity factors and the T-stress), special techniques, such as path-independent integrals, are necessary~\cite{kimpaulino2003,courtingardin2005,bergerkarageorghis2007,morais2007,passieuxgravouil2011}.  

\subsection{Approach} In this paper, we study the singular stress states that exist in solids when discontinuities are present in the geometry and/or the mechanical properties of the material using the scaled boundary formulation. The scaled boundary formulation is then combined with the XFEM to study the fracture parameters, for example, the stress intensity factors (SIFs) and the T-stress. The salient feature is that the combined method does not require a priori knowledge of the asymptotic expansions of the displacements. Moreover, the stiffness of the region containing the crack tip can be computed directly, unlike, the method described in~\cite{menkbordas2010,mousavigrinspun2011,duartekim2008,zhu2012}, which requires a two step procedure.

\subsection{Outline} The paper commences with a brief introduction to the extended finite element method (XFEM) in Section \ref{xfemoverview} and the scaled boundary finite element method (SBFEM) in Section \ref{sbfem}. The coupling of the XFEM and the SBFEM is discussed in Section \ref{couplexfemsbfem}. Section \ref{gsiftress} discusses the computation of the stress intensity factors (SIFs) and the T-stress in an isotropic and in a heterogeneous medium. Section \ref{numresults} presents the displacement and stress modes and the order of singularity for various configuration using the scaled boundary formulation. The effectiveness of the coupled scaled boundary and XFEM is demonstrated with a few problems from linear elastic fracture mechanics. Also, the formulation is used to model crack growth along the interface and into the material. The concluding remarks are given in the last section.

\section{Overview of the extended finite element method}

\label{xfemoverview}

\begin{figure}[htpb]
\centering
\scalebox{0.7}{\input{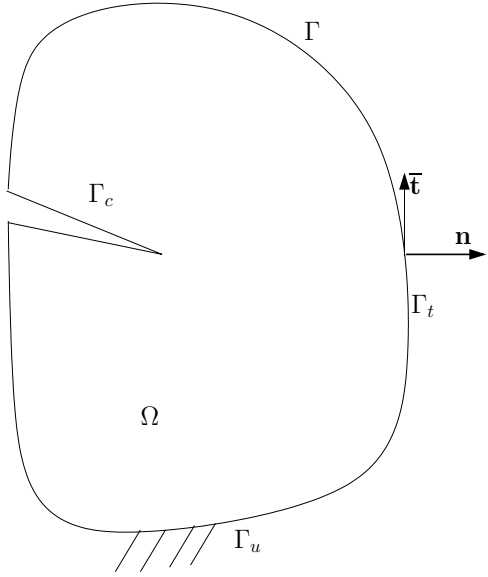}}
\caption{Two-dimensional elastic body with a crack.}
\label{fig:bodywcrack}
\end{figure}

Consider $\Omega \subset \mathbb{R}^2$, the reference configuration of a cracked linearly isotropic elastic body (see \fref{fig:bodywcrack}). The boundary of $\Omega$, denoted by $\Gamma$, is partitioned into three parts $\Gamma_u, \Gamma_n$ and $\Gamma_c$, where Dirichlet condition is prescribed on $\Gamma_u$ and Neumann condition is prescribed on $\Gamma_n$ and $\Gamma_c$.  The governing equilibrium equations for a 2D elasticity problem with internal boundary  $\Gamma_c$ defined in the domain $\Omega$ and bounded by $\Gamma$ is
\begin{equation}
\nabla_s^{\rm T} \bvsig + \bb = \mathbf{0} ~\hspace{0.5cm} \in ~\hspace{0.5cm} \Omega
\label{eqn:EE}
\end{equation}
where $\nabla_s(\cdot)$ is the symmetric part of the gradient operator, $\mathbf{0}$ is a null vector, $\bvsig$ is the stress tensor and $\bb$ is the body force. The boundary conditions for this problem are:
\begin{eqnarray}
\bvsig \cdot \cn = \overline{\mathbf{t}} \hspace{0.5cm} \textup{on} ~\hspace{0.5cm} \Gamma_t \nonumber \\
\uu = \overline{\uu} \hspace{0.5cm} \textup{on} ~\hspace{0.5cm} \Gamma_u \nonumber \\
\bvsig \cdot \cn = \overline{\mathbf{t}} \hspace{0.5cm} \textup{on} ~\hspace{0.5cm}  \Gamma_c
\end{eqnarray}
where $\overline{\uu} = (\overline{u}_x,\overline{u}_y)^{\rm T}$ is the prescribed displacement vector on the essential boundary $\Gamma_u$; $\overline{\mathbf{t}} = (\overline{t}_x,\overline{t}_y)^{\rm T}$ is the prescribed traction vector on the natural boundary $\Gamma_t$ and $\cn$ is the outward normal vector. In this study, it is assumed that the displacements remain small and the strain-displacement relation is given by $\bveps = \nabla_s \uu$. Let us assume a linear elastic behaviour, the constitutive relation is given by $\bvsig = \dd \colon \bveps$, where $\dd$ is a fourth order elasticity tensor. The weak formulation of the static problem is then given by:
\begin{equation}
\textup{find} ~\hspace{0.3cm} \uu^h \in \mathcal{U}^h \hspace{0.3cm} \textup{such~~that} \hspace{0.3cm} \forall \vv^h \in \mathcal{V}^h \hspace{0.3cm} a(\uu^h,\vv^h) = \ell(\vv^h)
\label{eqn:variweakform1}
\end{equation}
where $\mathcal{U}^h \subset \mathcal{U}$ and $\mathcal{V}^h \subset \mathcal{V}$.

\subsection{eXtended Finite Element Method (XFEM)}
\begin{figure}[htpb]
\centering
\includegraphics[angle=0,width=0.4\columnwidth]{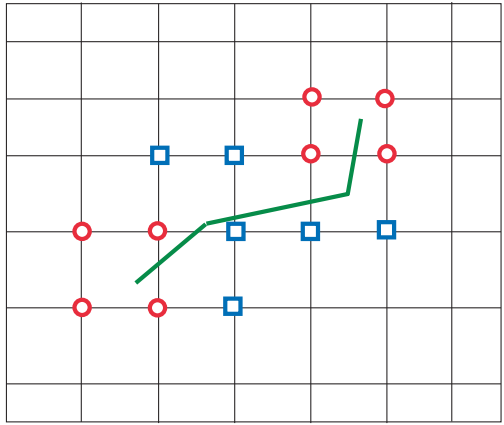}
\caption{A typical FE mesh with an internal discontinuity. Nodes are enriched with a Heaviside functions whose nodal support is cut by the crack interior (`\emph{squared}' nodes) and with set of functions that span the asymptotic fields whose nodal support is cut by the crack tip (`\emph{circled}' nodes).} 
\label{fig:xfemelementcategory}
\end{figure}

In the XFEM, the classical FEM polynomial approximation space is augmented by a set of functions called the `enrichment functions'. This is achieved through the framework of the partition of unity. XFEM is classified as one of the partition of unity methods (PUMs). A partition of unity in a domain $\Omega$ is a set of functions $N_I$ such that
\begin{equation}
\sum\limits_{I \in \mathcal{N}^{\rm {fem}}} N_I(\xx) = 1, \hspace{0.5cm} \xx \in \Omega
\end{equation}
where $\mathcal{N}^{\rm{fem}}$ is the set of nodes in the FE mesh. Using this property, any function $\varphi$ can be reproduced by a product of the PU shape functions with $\varphi$. Let $\uu^h \subset \mathcal{U}$, the XFEM approximation can be decomposed into the standard part $\uu^h_{\rm{fem}}$ and an enriched part $\uu^h_{\rm{enr}}$ as:
\begin{eqnarray}
\uu^h(\xx) &=& \uu^h_{\rm{fem}}(\xx) + \textup{Enrichment functions}
\label{eqn:xfemgeneralform}
\end{eqnarray}
where $\mathcal{N}^{\rm{fem}}$ is the set of all the nodes in the FE mesh. The modification of the displacement approximation does not introduce a new form of the discretized finite element equilibrium equations, but leads to an enlarged problem to solve. A generic form for the displacement approximation in case of the linear elastic fracture mechanics (LEFM) takes the form~\cite{belytschkoblack1999}:
\begin{equation}
\uu^h(\xx)=\sum\limits_{I \in \mathcal{N}^{\rm{fem}}} N_I(\xx)\qq_I + \sum\limits_{J \in \mathcal{N}^{\rm c}} N_J(\xx)\vartheta(\xx)\aaa_J + \sum\limits_{K \in
\mathcal{N}^{\rm f}} N_K (\xx)\sum_{\alpha=1}^n B_{\alpha}(r,\theta)\bb^{\alpha}_K ,\label{eqn:uS2}
\end{equation}
where $\mathcal{N}^{\rm c}$ is the set of nodes whose shape function support is cut by the crack interior (`squared' nodes in \fref{fig:xfemelementcategory}) and $\mathcal{N}^{\rm f}$ is the set of nodes whose shape function support is cut by the crack tip (`circled' nodes in \fref{fig:xfemelementcategory}). $\vartheta$ and $B_{\alpha}$ are the enrichment functions chosen to capture the displacement jump across the crack surface and the singularity at the crack tip and $n$ is the total number of asymptotic functions. $\aaa_J$ and $\bb^{\alpha}_K$ are the nodal degrees of freedom corresponding to functions $\vartheta$ and $B_{\alpha}$, respectively. $n$ is the total number of near-tip asymptotic functions and $(r,\theta)$ are the local crack tip coordinates. A priori knowledge of the shape functions is required for the successful application of the method. Since its inception, this has attracted researchers to investigate the enrichment functions that can accurately describe the asymptotic fields ahead of the crack tip~\cite{asadpouremohammadi2007,hattorirojas-diaz2012,asharimohammadi2011,bouhalashao2013}.



\section{Overview of the scaled boundary formulation}
\label{sbfem}
The scaled boundary finite element method (SBFEM) is a semi-analytical computational technique developed by Wolf and Song~\cite{wolfsong2001}. It reduces the governing partial differential equations to a set of ordinary differential equations and is suitable for solving linear elliptic, parabolic and hyperbolic partial differential equations. In the SBFEM, the coordinates are defined by a radial-circumferential like coordinate system. Numerical solutions are sought around the circumferential direction using the conventional FEM, whilst in the radial direction, the solution is defined by smooth analytical functions. A scaling centre $O$ is selected at a point from which the whole boundary of the domain is visible (see \fref{fig:sbfemdescrip}). The radial coordinate $\xi$ is defined from the scaling centre $(\xi=0)$ to the boundary $(\xi=1)$. Each edge on the boundary is discretised using 1D finite elements with local coordinate $\eta$, defined in the range $-1\leq\eta\leq1$. The coordinate of a point on the element $\mathbf{x}_{\eta}(\eta)=[\begin{array}{cc} x(\eta) & y(\eta)\end{array}]^{\mathrm{T}}$ are expressed as 
\begin{equation}
\mathbf{x}_{\eta}(\eta)=\mathbf{N}(\eta)\mathbf{x}_{\mathrm{b}}\label{eq:lineEleGeom}
\end{equation}
where $\xx_b$ are the coordinates of the nodes on the boundary and $\mathbf{N}(\eta)$ is the shape function matrix. In this study, standard 1D Gauss-Lobatto-Lagrange shape functions are used. By choosing the origin of the Cartesian coordinate system as the scaling center, the geometry of the domain is described as:
\begin{equation}
\mathbf{x}(\xi,\eta)=\xi\mathbf{x}_{\eta}(\eta)\label{eqn:sbccartesian}
\end{equation}
\begin{figure}[htpb]
\centering
\subfigure[]{\scalebox{0.6}{\input{./Figures/sbfemnocrk.pstex_t}}}
\subfigure[]{\scalebox{0.6}{\input{./Figures/sbfemcrk.pstex_t}}}
\caption{Coordinate transformation in the scaled boundary finite element method: (a) domain without a crack. The scaling centre $O$ can be placed anywhere as long as it satisfies the visibility criteria and (b) domain with a crack and a material interface. The scaling center $O$ is placed at the crack tip. The material interface within the domain need not be discretized, while on the boundary a conforming mesh is built with a node at the material interface.}
\label{fig:sbfemdescrip}
\end{figure}
Using \Eref{eqn:sbccartesian}, the polar coordinates $r$ and $\theta$ are expressed in the scaled boundary coordinates $\xi$ and $\eta$ as 
\begin{align}
r(\xi,\eta) & =\xi r_{\eta}(\eta)=\xi\sqrt{x_{\eta}^{2}(\eta)+y_{\eta}^{2}(\eta)}\label{eq-polar-r}\\
\theta(\eta) & =\arctan\dfrac{y_{\eta}(\eta)}{x_{\eta}(\eta)}\label{eq-polar-theta}
\end{align}
where $r_{\eta}(\eta)$ is the radial coordinate on the boundary. The angle $\theta(\eta)$ is a single-valued function in its principal value ($-\pi<\theta\leq\pi$). The element number and the local coordinate $\eta$ can be regarded as a discrete representation of the angle $\theta$. 

\paragraph*{Displacement approximation} The displacements of a point (for example, point $A$ see \fref{fig:sbfemdescrip}) is approximated by: 
\begin{equation}
\mathbf{u}(\xi,\eta)=\mathbf{N}(\eta)\mathbf{u}(\xi)\label{eqn:dispapprox}
\end{equation}
 where $\mathbf{u}(\xi)$ are the radial displacement functions. Substituting \Eref{eqn:dispapprox} in the definition of strain-displacement relations, the strains $\boldsymbol{\varepsilon}(\xi,\eta)$ are expressed as:
\begin{equation}
\boldsymbol{\varepsilon}(\xi,\eta)=\mathbf{L}\mathbf{u}(\xi,\eta)\label{eqn:sbfemstrain}
\end{equation}
 where $\mathbf{L}$ is a linear operator matrix formulated in the
scaled boundary coordinates as 
\begin{equation}
\mathbf{L}=\mathbf{b}_{1}(\eta)\frac{\partial}{\partial\xi}+\xi^{-1}\mathbf{b}_{2}(\eta)\label{eqn:Loperator}
\end{equation}
 with 
\begin{align}
\mathbf{b}_{1}(\eta) & =\frac{1}{|\mathbf{J}(\eta)|}\left[\begin{array}{cc}
y_{\eta}(\eta)_{,\eta} & 0\\
0 & -x_{\eta}(\eta)_{,\eta}\\
-x_{\eta}(\eta)_{,\eta} & y_{\eta}(\eta)_{,\eta}
\end{array}\right] \nonumber \\
\mathbf{b}_{2}(\eta) & =\frac{1}{|\mathbf{J}(\eta)|}\left[\begin{array}{cc}
-y_{\eta}(\eta) & 0\\
0 & x_{\eta}(\eta)\\
x_{\eta}(\eta) & y_{\eta}(\eta)
\end{array}\right]\label{eq:b2}
\end{align}
 In \Eref{eq:b2}, $\mathbf{J}(\eta)$ is the Jacobian on the boundary \cite{wolfsong2001,deekswolf2002}, given by:
\begin{align*}
\mathbf{J}(\eta)= & \left[\begin{array}{cc}
x_{\eta}(\eta) & y_{\eta}(\eta)\\
x_{\eta}(\eta)_{,\eta} & y_{\eta}(\eta)_{,\eta}
\end{array}\right]
\end{align*}
Using ~\Eref{eqn:sbfemstrain} and the Hooke's law $\boldsymbol{\sigma}=\mathbf{D}\boldsymbol{\varepsilon}$, the stresses $\boldsymbol{\sigma}(\xi,\eta)$ is expressed as 
\begin{equation}
\boldsymbol{\sigma}(\xi,\eta)=\mathbf{D}\left(\mathbf{B}_{1}(\eta)\mathbf{u}(\xi)_{,\xi}+\xi^{-1}\mathbf{B}_{2}(\eta)\mathbf{u}(\xi)\right)\label{eqn:sbfemstress}
\end{equation}
 where $\mathbf{D}$ is the material constitutive matrix and 
\begin{align}
\mathbf{B}_{1}(\eta) & =\mathbf{b}_{1}(\eta)\mathbf{N}(\eta) \nonumber \\
\mathbf{B}_{2}(\eta) & =\mathbf{b}_{2}(\eta)\mathbf{N}(\eta)_{,\eta}
\end{align}
Substituting~\Erefs{eqn:dispapprox}, (\ref{eqn:sbfemstrain}) and (\ref{eqn:sbfemstress}) for the displacement, strain and stress fields, respectively, in the virtual work statement results in \cite{wolfsong2001,deekswolf2002}
\begin{equation}
\begin{split}\delta\mathbf{u}{}_{\mathrm{b}}^{{\rm T}}\left((\mathbf{E}_{0}\xi\mathbf{u}(\xi)_{,\xi}+\mathbf{E}_{1}^{{\rm T}}\mathbf{u}(\xi))|_{\xi=1}-\mathbf{f}\right)\\
-\int\limits _{0}^{1}\delta\mathbf{u}(\xi)^{{\rm T}}\left(\mathbf{E}_{0}\xi^{2}\mathbf{u}(\xi)_{,\xi\xi}+(\mathbf{E}_{0}+\mathbf{E}_{1}^{\mathrm{T}}-\mathbf{E}_{1})\xi\mathbf{u}(\xi)_{,\xi}-\mathbf{E}_{2}\mathbf{u}(\xi)\right)d\xi=0
\end{split}
\label{eqn:sbfemvirtual}
\end{equation}
where $\mathbf{f}$ is the equivalent boundary nodal forces and $\mathbf{u}_{\mathrm{b}}$ is the nodal displacement vector. By considering the arbitrariness of $\delta\mathbf{u}(\xi)$, the following ODE is obtained: 
\begin{equation}
\mathbf{E}_{0}\xi^{2}\mathbf{u}(\xi)_{,\xi\xi}+(\mathbf{E}_{0}+\mathbf{E}_{1}^{\mathrm{T}}-\mathbf{E}_{1})\xi\mathbf{u}(\xi)_{,\xi}-\mathbf{E}_{2}\mathbf{u}(\xi)=0\label{eqn:governODEsbfem}
\end{equation}
where $\mathbf{E}_{0},\mathbf{E}_{1}$ and $\mathbf{E}_{2}$ are coefficient matrices given by \cite{wolfsong2001,deekswolf2002} 
\begin{align}
\mathbf{E}_{0} & =\int_{\eta}\mathbf{B}_{1}(\eta)^{{\rm T}}\mathbf{D}\mathbf{B}_{1}(\eta)|\mathbf{J}(\eta)|d\eta,\nonumber \\
\mathbf{E}_{1} & =\int_{\eta}\mathbf{B}_{2}(\eta)^{{\rm T}}\mathbf{D}\mathbf{B}_{1}(\eta)|\mathbf{J}(\eta)|d\eta,\nonumber \\
\mathbf{E}_{2} & =\int_{\eta}\mathbf{B}_{2}(\eta)^{{\rm T}}\mathbf{D}\mathbf{B}_{2}(\eta)|\mathbf{J}(\eta)|d\eta.\label{eqn:coeffmat}
\end{align}
They are evaluated element-by-element on the element boundary and assembled similar to the standard FE procedure of assemblage. For a 2 noded line element, the coefficient matrices given by~\Eref{eqn:coeffmat} can be written explicitly. A simple routine is given in~\cite{natarajansong2013}. The code can easily be extended to treat higher order elements.

\paragraph*{Computation of the stiffness matrix}
\Eref{eqn:governODEsbfem} is a homogeneous second-order ordinary differential equation. Its solution is obtained by introducing the variable $\boldsymbol{\chi}(\xi)$
\begin{align}
\boldsymbol{\chi}(\xi)= & \left\{ \begin{array}{c}
\mathbf{u}(\xi)\\
\mathbf{q}(\xi)
\end{array}\right\} \label{eq:chi}
\end{align}
 where $\mathbf{q}(\xi)$ is the internal load vector 
\begin{align}
\mathbf{q}(\xi)= & \mathbf{E}_{0}\xi\mathbf{u}(\xi)_{,\xi}+\mathbf{E}_{1}^{{\rm T}}\mathbf{u}(\xi)
\end{align}
The boundary nodal forces are related to the displacement functions by: 
\begin{equation}
\mathbf{f}=\mathbf{q}(\xi=1)=(\mathbf{E}_{0}\xi\mathbf{u}(\xi)_{,\xi}+\mathbf{E}_{1}^{{\rm T}}\mathbf{u}(\xi))|_{\xi=1}\label{eqn:nodalforce}
\end{equation}
This allows~\Eref{eqn:governODEsbfem} to be transformed into a first order ordinary differential equation with twice the number of unknowns in an element as: 
\begin{equation}
\xi\boldsymbol{\chi}(\xi)_{,\xi}=-\mathbf{Z}\boldsymbol{\chi}(\xi)\label{eq:first order}
\end{equation}
 where $\mathbf{Z}$ is a Hamiltonian matrix 
\begin{equation}
\mathbf{Z}=\left[\begin{array}{cc}
\mathbf{E}_{0}^{-1}\mathbf{E}_{1}^{\mathrm{T}} & -\mathbf{E}_{0}^{-1}\\
\mathbf{E}_{1}\mathbf{E}_{0}^{-1}\mathbf{E}_{1}^{\mathrm{T}}-\mathbf{E}_{2} & -\mathbf{E}_{1}\mathbf{E}_{0}^{-1}
\end{array}\right]\label{eq:Hamiltonian matrix}
\end{equation}
An eigenvalue decomposition of $\mathbf{Z}$ is performed. The blocks of eigenvalues and transformation matrices necessary are:
\begin{align}
\mathbf{Z}\left[\begin{array}{c}
\boldsymbol{\Phi}_{\mathrm{u}}\\
\boldsymbol{\Phi}_{\mathrm{q}}
\end{array}\right] & =\left[\begin{array}{c}
\boldsymbol{\Phi}_{\mathrm{u}}\\
\boldsymbol{\Phi}_{\mathrm{q}}
\end{array}\right]\boldsymbol{\Lambda}_{\mathrm{n}}\label{eq:eigen decomp}
\end{align}
In \Eref{eq:eigen decomp}, $\boldsymbol{\Lambda}_{\mathrm{n}}=\mathrm{diag}\left(\lambda_{1},\,\lambda_{2},\,...,\lambda_{n}\right)$ contains the eigenvalues with negative real part. $\boldsymbol{\Phi}_{\mathrm{u}}$ and $\boldsymbol{\Phi}_{\mathrm{q}}$ are the corresponding transformation matrices of $\boldsymbol{\Lambda}_{\mathrm{n}}$. They represent the modal displacements and forces, respectively. The general solution of \Eref{eq:first order} is given by:
\begin{align}
\mathbf{u}(\xi)= & \boldsymbol{\Phi}_{\mathrm{u}}\xi^{-\boldsymbol{\Lambda}_{\mathrm{n}}}\mathbf{c}\\
\mathbf{q}(\xi)= & \boldsymbol{\Phi}_{\mathrm{q}}\xi^{-\boldsymbol{\Lambda}_{\mathrm{n}}}\mathbf{c}
\label{eqn:raddispfunc}
\end{align}
where $\mathbf{c}$ are integration constants that are obtained from the nodal displacements $\mathbf{u}_{\mathrm{b}}=\mathbf{u}(\xi=1)$ as:
\begin{align}
\mathbf{c}= & \boldsymbol{\Phi}_{\mathrm{u}}^{-1}\mathbf{u}_{\mathrm{b}}\label{eq:int constants}
\end{align}
The complete displacement field of a point defined by the sector covered by a line element on the element is obtained by substituting~\Eref{eqn:raddispfunc} into \Eref{eqn:dispapprox} resulting in:
\begin{align}
\mathbf{u}(\xi,\eta)= & \mathbf{N}(\eta)\boldsymbol{\Phi}_{\mathrm{u}}\xi^{-\boldsymbol{\Lambda}_{\mathrm{n}}}\mathbf{c}\label{eq:dispapprox final}
\end{align}
Taking the derivative of $\mathbf{u}(\xi)$ with respect to $\xi$ and substituting into~\Eref{eqn:sbfemstress} the stress field $\boldsymbol{\sigma}(\xi,\eta)$ can be expressed as:
\begin{align}
\boldsymbol{\sigma}(\xi,\eta)= & \boldsymbol{\Psi}_{\sigma}(\eta)\xi^{-\boldsymbol{\Lambda}_{\mathrm{n}}-\mathbf{I}}\mathbf{c}\label{eq:stress field complete}
\end{align}
where the stress mode $\boldsymbol{\Psi}_{\sigma}(\eta)$ is defined as: 
\begin{align}
\boldsymbol{\Psi}_{\sigma}(\eta)= & \mathbf{D}\left(-\mathbf{B}_{1}(\eta)\boldsymbol{\Phi}_{\mathrm{u}}\boldsymbol{\Lambda}_{\mathrm{n}}+\mathbf{B}_{2}(\eta)\boldsymbol{\Phi}_{\mathrm{u}}\right)\label{eq:stress mode}
\end{align}
The stiffness matrix of an element is obtained by first substituting \Eref{eq:int constants} into~\Eref{eqn:raddispfunc} at $\xi=1$. This results in:
\begin{align}
\mathbf{f}= & \boldsymbol{\Phi}_{\mathrm{q}}\boldsymbol{\Phi}_{\mathrm{u}}^{-1}\mathbf{\mathbf{u}_{\mathrm{b}}}\label{eq:equil poly}
\end{align}
From \Eref{eq:equil poly}, the stiffness matrix $\mathbf{K}$ can be identified to be given by the expression
\begin{equation}
\mathbf{K}=\boldsymbol{\Phi}_{\mathrm{q}}\boldsymbol{\Phi}_{\mathrm{u}}^{-1}\label{eqn:sbfemkmat-b}
\end{equation}

\begin{remark}
No a priori knowledge of the asymptotic solution and special integration technique is required to compute the stiffness matrix. The SBFEM solution includes the stress singularity as well. The details are given in Section \ref{gsiftress}.
\end{remark}

\begin{remark}
The SBFEM can be formulated over arbitrary regions, for example over polygonal elements~\cite{ooisong2012,natarajanooi2013}.
\end{remark}

\section{Coupling of the XFEM \& the SBFEM} \label{couplexfemsbfem}
\begin{figure}[htpb]
\centering
\includegraphics[scale=0.5]{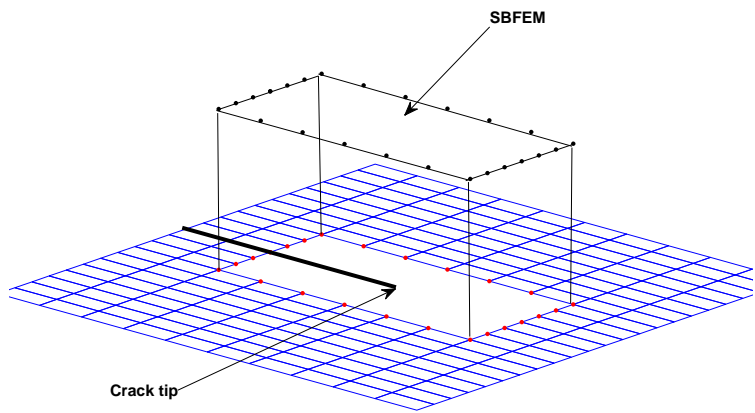}
\caption{Coupling of XFEM and SBFEM}
\label{fig:couplexfemsbfem}
\end{figure}

Recently, the authors proposed the extended SBFEM (xSBFEM)~\cite{natarajansong2013}. The main idea of the xSBFEM is to approximate the non-smooth behaviour in the close vicinity by the scaled boundary formulation instead of augmenting the approximation basis with asymptotic expansion. The nodes whose nodal support is completely cut by the crack interior or the material interface are treated within the XFEM framework and enriched with a Heaviside function or 'abs' enrichment function. The scaling center is chosen as the crack tip (see \fref{fig:couplexfemsbfem}). This local modification and enrichment strategy introduces four types of elements:

\begin{itemize}
\item \textit{Split elements} are elements completely cut by the crack interior or by the material interface. Their nodes are enriched with the Heaviside function to represent the jump across the crack face or by the 'abs' function to represent the weak discontinuity.
\item \textit{Split-blending elements} are elements neighbouring split elements. They are such that some of their nodes are enriched with the weakly discontinuous function and others are not enriched at all.
\item \textit{Scaled boundary elements} are elements in the close proximity of the crack tip treated by the semi-analytical formulation.
\item \textit{Standard elements} are elements that are in neither of the above categories. None of their nodes are enriched.
\end{itemize}

\begin{remark}
It is noted that no special treatment is required for the elements that share nodes between the $\Omega^{\rm fem}$ and $\Omega^{^{\rm sbfem}}$.
\end{remark}

The displacement approximation can be decomposed into the standard part $\uu^h_{\rm std}$ and an enriched part $\uu^h_{\rm xfem}$ as given by \Eref{eqn:uS2} without asymptotic enrichments. In the proposed approach, the standard part is further decomposed into the standard FEM $\uu^h_{\rm fem}$ and into a semi-analytical part $\uu^h_{\rm sbfem}$. The displacement approximation is then given by:
\begin{equation}
\renewcommand{\arraystretch}{2}
\uu^h(\xx) = \left\{ \begin{array}{lr} \sum\limits_{I \in \mathcal{N}^{\rm fem}} N_I(\xx) \qq_I + \sum\limits_{J \in \mathcal{N}^c} N_J(\xx) H(\xx) \aaa_J & \xx \in \Omega^{^{\rm fem}} \cup \Omega^{^{\rm xfem}} \\
 \sum\limits_{K \in \mathcal{N}^{\rm sbfem}} N_K(\eta)  \uu_{bK}(\xi) & \xx = f(\xi,\eta) \in \Omega^{^{\rm sbfem}} \end{array} \right.
 \label{eqn:dispApprox}
\end{equation}
where $\mathcal{N}^{\rm fem}$ is the set of all nodes in the FE mesh that belongs to $\Omega^{^{\rm fem}}$, $\mathcal{N}^{\rm sbfem}$ is the set of nodes in the FE mesh that belongs to $\Omega^{^{\rm sbfem}}$ and $\mathcal{N}^{\rm xfem}$ is the set of nodes that are enriched with the Heaviside function $H$. $N_I$ and $N_J$ are the standard finite element functions, $\qq_I$ and $\aaa_J$ are the standard and the enriched nodal variables associated with node $I$ and node $J$, respectively. $\uu_{bK}$ are the nodal variables associated with the nodes on the boundary of $\Omega^{^{\rm sbfem}}$ and $N_K$ are the standard FE shape functions defined on the scaled boundary coordinates in $\Omega^{^{\rm sbfem}}$. For the nodes that on the boundary of $\Omega^{^{\rm fem}}$ and $\Omega^{^{\rm sbfem}}$, no special coupling technique is required. As the nodal unknown coefficients are required to be continuous across the boundary, the unknown coefficients from the SBFEM and the FEM are the same and are assembled in the usual way. A similar procedure is followed when assembling the $\Omega^{^{\rm fem}}$ and $\Omega^{^{\rm xfem}}$ to the global stiffness and to the global force vector. For more details, interested readers are referred to~\cite{natarajansong2013}.

\section{Computation of the stress intensity factors and the T-Stress} \label{gsiftress}
An attractive feature of the SBFEM is that no a priori knowledge of the asymptotic solution is required to accurately handle the stress singularity at a crack tip. When the scaling centre is placed at the crack tip, the solution for the stress field in \Eref{eq:stress field complete} is expressed, by using~\Eref{eq-polar-r}, as

\begin{equation}
\boldsymbol{\sigma}(r,\eta)=\sum_{i=1}^{n}c_{i}r^{-(\lambda_{i}+1)}\left(r_{\eta}^{\lambda_{i}+1}(\eta)\boldsymbol{\psi}_{\sigma i}(\eta)\right)\label{eq-stress-bd-slu}
\end{equation}
where $\boldsymbol{\psi}_{\sigma i}(\eta)$  is the $i-$th stress mode, i.e. the $i-$th column of the matrix $\boldsymbol{\Psi}_{\sigma}(\eta)$. Like the well-known William expansion~\cite{williams1957}, \Eref{eq-stress-bd-slu} is a power series of the radial coordinate $r$. The radial variation of each term of the series is expressed analytically by the power function $r^{-(\lambda_{i}+1)}$. At discrete points along the boundary, the angular coordinates (see~\Eref{eq-polar-theta}) are arranged as a vector $\boldsymbol{\theta}(\eta)$ and the stress modes $\boldsymbol{\psi}_{\sigma i}(\eta)$ are computed. $\boldsymbol{\psi}_{\sigma i}(\eta)$ and $\boldsymbol{\theta}(\eta)$ form a parametric equation of the angular variation of stresses. The singular stress and the T-stress terms can be easily identified by the value of the exponent $-(\lambda_{i}+1)$. When the real part of the exponent $-(\lambda_{i}+1)$ of a term is negative, the stresses of this term at the crack tip, i.e. $\xi=0$, tend to infinity. When the exponent $-(\lambda_{i}+1)$ of a term is equal to $0$, the stresses of this term are constant and contribute to the T-stress. 

\begin{figure}[htpb]
\centering
\scalebox{0.8}{\input{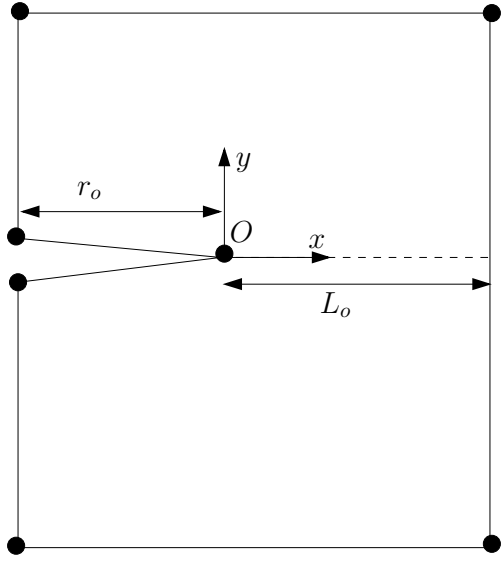}}
\caption{A cracked domain modelled by the scaled boundary finite element method. $L_{o}=r_{\eta}(\theta=0)$ is the distance from the crack tip to the boundary. At the boundary $\xi=$ 1.}
\label{fig:crkpolysbfem}
\end{figure}

\paragraph{Crack in homogeneous material or material interface} In the case of a crack in a homogeneous material or on a material interface, two singular terms exist in the solution. Denoting the singular stress modes as $I$ and $II$, the singular stresses $\bvsig^s(\xi,\eta)$ (superscript $\mathrm{s}$ for singular stresses) are obtained from \Eref{eq-stress-bd-slu}:
\begin{equation}
\boldsymbol{\sigma}^{\mathrm{s}}(r,\eta)=\sum_{i=I,II}c_{i}r^{-(\lambda_{i}+1)}\left(r_{\eta}^{\lambda_{i}+1}(\eta)\boldsymbol{\psi}_{\sigma i}(\eta)\right)\label{eq:singularStress}
\end{equation}
Note that the singular stress terms are separated from other terms and the stress singularity is represented analytically. This allows the evaluation of the stress intensity factors by directly matching their definition with the singular stresses. For convenience, the point on the boundary along the crack front $\theta=0$ (see \fref{fig:crkpolysbfem}) is considered. The distance from the crack tip to the boundary is denoted as $L_{o}=r_{\eta}(\theta=0)$. From \Eref{eq:singularStress}, the values of the singular stresses at this point are equal to: 
\begin{equation}
\boldsymbol{\sigma}^{\mathrm{s}}(L_{o},\theta=0)=\sum_{i=I,II}c_{i}\boldsymbol{\psi}_{\sigma i}(\theta=0)\label{eq:singularStress_theta0}
\end{equation}
where $\boldsymbol{\psi}_{\sigma i}(\theta=0)$ is the value of the stress modes at $\theta=0$. It is obtained by interpolating $\boldsymbol{\psi}_{\sigma i}(\eta)$ at the discrete points of $\boldsymbol{\theta}(\eta)$.

Stress intensity factors can be computed directly from their definition using the stresses in ~\Eref{eq:singularStress_theta0}. For a crack in a homogeneous medium, the classical definition of stress intensity factors $K_{I}$ and $K_{II}$ for mode I and II are expressed as:
\begin{equation}
\left\{ \begin{array}{c}
K_{I}\\
K_{II}
\end{array}\right\} =\sqrt{2\pi r}\left\{ \begin{array}{c}
\sigma_{\theta\theta}^{\mathrm{s}}(r,\theta=0)\\
\tau_{r\theta}^{\mathrm{s}}(r,\theta=0)
\end{array}\right\} \label{eq:GSIFdef-case1-2}
\end{equation}
Formulating ~\Eref{eq:GSIFdef-case1-2} at $r=L_{o}$ results in
\begin{equation}
\left\{ \begin{array}{c}
K_{I}\\
K_{II}
\end{array}\right\} =\sqrt{2\pi L_{0}}\left\{ \begin{array}{c}
\sigma_{\theta\theta}^{\mathrm{s}}(L_{o},\theta=0)\\
\tau_{r\theta}^{\mathrm{s}}(L_{o},\theta=0)
\end{array}\right\} \label{eq:GSIFdef-case1-2-1}
\end{equation}
The stress intensity factors are then determined by substituting the stress components $\sigma_{\theta\theta}^{\mathrm{s}}(L_{o},\theta=0)$ and $\tau_{r\theta}^{\mathrm{s}}(L_{o},\theta=0)$ obtained from \Eref{eq:singularStress_theta0} into \Eref{eq:GSIFdef-case1-2-1}. 

\paragraph{Crack at an interface} In the case of a crack at the interface between two isotropic materials, the orders of singularity are a pair of complex conjugates $0.5\pm\mathrm{i}\epsilon$, where the oscillatory index $\epsilon$ depends on the material properties. The stress intensity factors are defined as:

\begin{equation}
\left\{ \begin{array}{c}
K_{I}\\
K_{II}
\end{array}\right\} =\sqrt{2\pi r}\left[\begin{array}{cc}
\cos(\epsilon\ln(r/L)) & \sin(\epsilon\ln(r/L))\\
-\sin(\epsilon\ln(r/L)) & \cos(\epsilon\ln(r/L))
\end{array}\right]\left\{ \begin{array}{c}
\sigma_{\theta\theta}^{\mathrm{(s)}}(r,0)\\
\tau_{r\theta}^{\mathrm{(s)}}(r,0)
\end{array}\right\} \label{eq:GSIFdef-case2-3}
\end{equation}
where $L$ is the characteristic length. Formulating the definition (see~\Eref{eq:GSIFdef-case2-3}) at $r=L_{o}$, the stress intensity factors are expressed as
\begin{equation}
\left\{ \begin{array}{c}
K_{I}\\
K_{II}
\end{array}\right\} =\sqrt{2\pi L_{o}}\left[\begin{array}{cc}
\cos(\epsilon\ln(L_{o}/L)) & \sin(\epsilon\ln(L_{o}/L))\\
-\sin(\epsilon\ln(L_{o}/L)) & \cos(\epsilon\ln(L_{o}/L))
\end{array}\right]\left\{ \begin{array}{c}
\sigma_{\theta\theta}^{\mathrm{(s)}}(L_{o},0)\\
\tau_{r\theta}^{\mathrm{(s)}}(L_{o},0)
\end{array}\right\} \label{eq:GSIFdef-case2-3-1}
\end{equation}
where the singular stresses are obtained from \Eref{eq:singularStress_theta0}. Although, it is out of the scope of this paper, it is worthwhile to mention that the scaled boundary finite element method is capable of modeling all types of stress singularity occurring at multi-material junctions~\cite{song2005}.

For a subdomain containing a crack tip, two of the eigenvalues are equal to $1$. They represent the $T-$stress term and the rotational rigid body motion term, which does not contribute to the stresses. They are separated from other terms in \Eref{eq:singularStress_theta0} and expressed as (superscript T for the $T-$stress) 
\begin{equation}
\boldsymbol{\sigma}^{T}(\eta)=\sum_{i=T_{I},T_{II}}c_{i}\boldsymbol{\psi}_{\sigma i}(\eta)\label{eq:TStress-1}
\end{equation}
The T-stress along the crack front ($\theta=0$) is determined by interpolating the angular variation of the two stress modes $(\boldsymbol{\theta}(\eta),\boldsymbol{\psi}_{\sigma i}(\eta))$. 

\begin{remark}
Note that when evaluating the stress intensity factors or the $T-$stress, only the corresponding modes are involved. It is not necessary to evaluate the complete stress field. 
\end{remark}

\section{Results and discussion} 
\label{numresults}
In this section, we first employ the SBFEM to estimate the order of the singularity and the stress and the displacement modes for known benchmark problems in the context of linear elastic fracture mechanics. Later, we will adopt the combined XFEM and the SBFEM to compute the stress intensity factors and the T-stress for cracks aligned to the material interfaces. The stress distribution ahead of the crack and the T-stress for a crack terminating at an interface is also studied. The results of which are compared with a finite element solution. In this study, unless otherwise mentioned, bilinear quadrilateral elements are used. An extension to other arbitrary polytopes is possible as illustrated in~\cite{natarajanooi2013}. The results from the present approach are compared with available numerical results in the literature, for example, with the XFEM~\cite{yuwu2012} and with the BEM~\cite{sladeksladek1997,kimvlassak2006}.

\subsection{Computation of displacement \& stress modes and the order of singularity}
To test and to demonstrate the convergence and the accuracy of the SBFEM in estimating the strength of the singularity and the corresponding displacement and stress modes, we consider three examples. 

\begin{figure}[htpb]
\centering
\subfigure[]{\scalebox{0.5}{\input{./Figures/isocrk.pstex_t}}}
\subfigure[]{\includegraphics[scale=0.5]{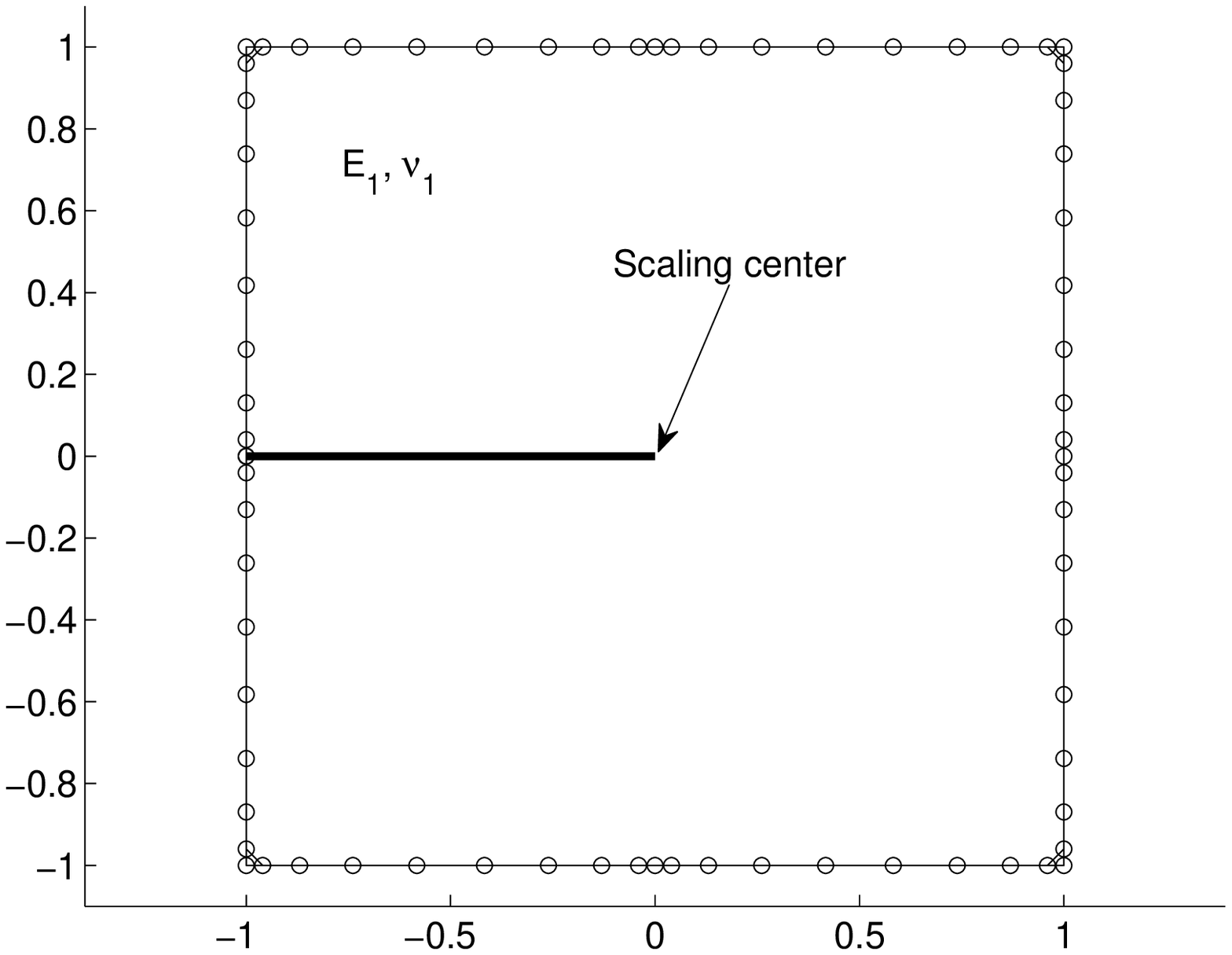}}
\caption{Crack in an isotropic medium: (a) geometry and (b) SBFEM Mesh detail}
\label{fig:isoprobdescription}
\end{figure}

\paragraph{Crack in an isotropic medium} \fref{fig:isoprobdescription} shows a crack in an isotropic medium with Young's modulus $E$ and Poisson's ratio $\nu$. The displacement fields are given by~\cite{williams1957}:
\begin{align}
u_x(r,\theta) &= \frac{K_I}{2\mu} \sqrt{ \frac{r}{2\pi}} \cos \frac{\theta}{2} \left[ \kappa -1 + 2 \sin^2\frac{\theta}{2} \right] + \frac{K_{II}}{2\mu} \sqrt{ \frac{r}{2\pi}} \sin \frac{\theta}{2} \left[ \kappa + 1 + 2\cos^2 \frac{\theta}{2} \right] \nonumber \\
u_y(r,\theta) &= \frac{K_I}{2\mu} \sqrt{ \frac{r}{2\pi}} \sin \frac{\theta}{2} \left[ \kappa +1 - 2 \cos^2\frac{\theta}{2} \right] + \frac{K_{II}}{2\mu} \sqrt{ \frac{r}{2\pi}} \sin \frac{\theta}{2} \left[ \kappa - 1 - 2\sin^2 \frac{\theta}{2} \right]
\end{align}
where $\kappa$ is the Kolosov constant and $\mu$ is the shear modulus. The reference solution is obtained along a circle of unit radius $r=$ 1 and with stress intensity factors $K_I$ is set to unity and $K_{II}$ to zero. The displacement modes, the stress modes and the order of singularity are computed by the scaled boundary formulation. A typical scaled boundary finite element mesh is shown in \fref{fig:isoprobdescription}. It is noted that only the boundary of the domain needs to be discretized and the solution within the domain is represented analytically. 

\begin{remark}
In our study, we have employed spectral element to discretize the boundary. We employ Lagrange interpolants within each element, where the nodes of these shape functions are placed at the zeros of Legendre polynomials (Gauss Lobatto points) mapped from the reference domain $[-1,1]\times[-1,1]$ to each element.
\end{remark}

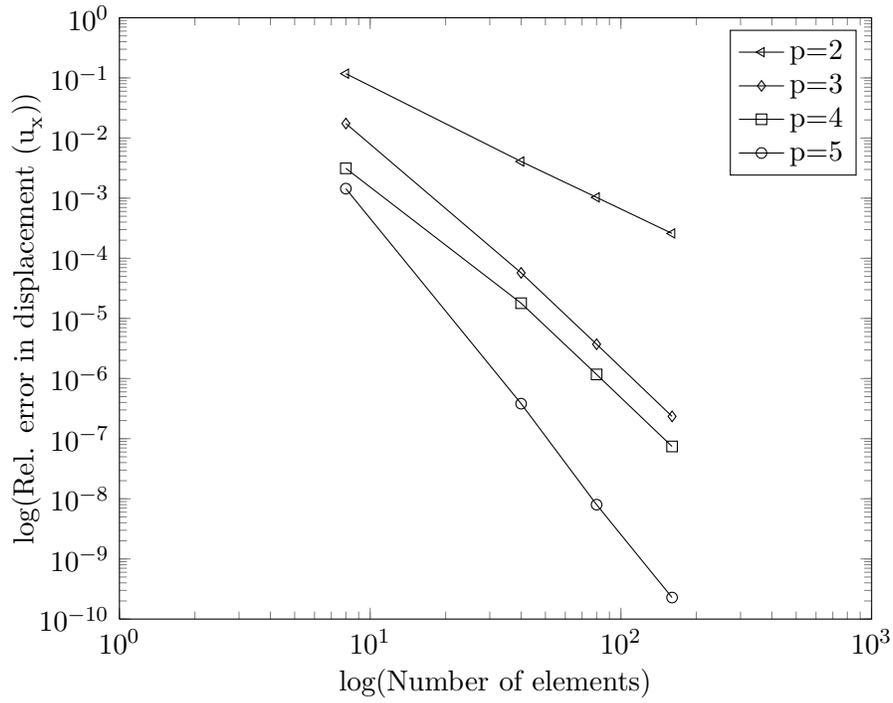
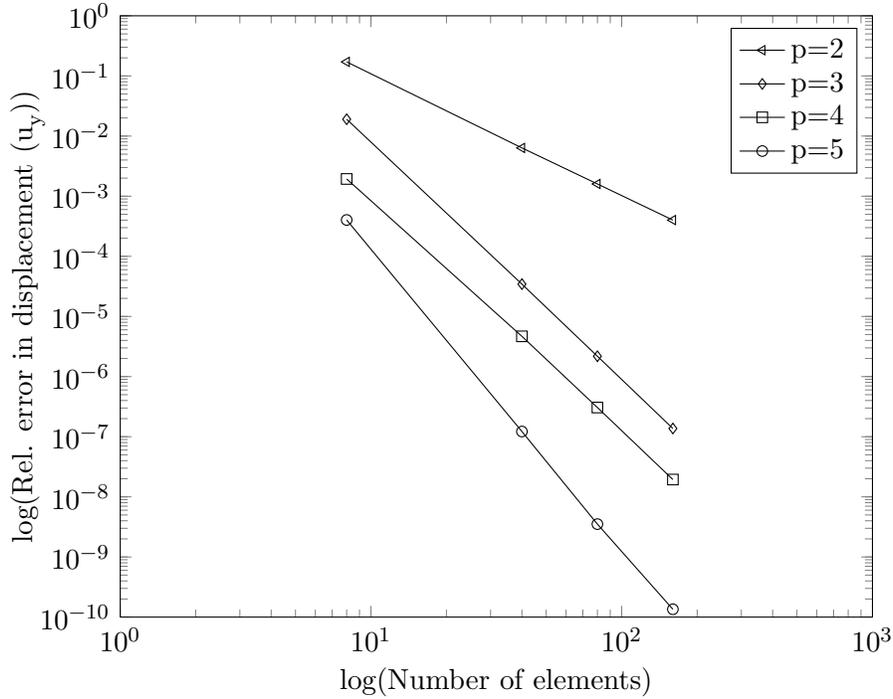
\begin{figure}
\centering
\newlength\figureheight 
\newlength\figurewidth 
\setlength\figureheight{8cm} 
\setlength\figurewidth{10cm}
\subfigure[Relative error in displacement, $u_x$]{
%
%
%
%
\begin{tikzpicture}

\begin{axis}[%
width=\figurewidth,
height=\figureheight,
scale only axis,
xmode=log,
xmin=1,
xmax=1000,
xminorticks=true,
xlabel={log(Number of elements)},
ymode=log,
ymin=1e-10,
ymax=1,
yminorticks=true,
ylabel={$\text{log(Rel. error in displacement (u}_\text{x}\text{))}$},
legend style={draw=black,fill=white,legend cell align=left}
]
\addplot [
color=black,
solid,
mark=triangle,
mark options={solid,,rotate=90}
]
table[row sep=crcr]{
8 0.117541408\\
40 0.004070536\\
80 0.001029892\\
160 0.000258761\\
};
\addlegendentry{p=2};

\addplot [
color=black,
solid,
mark=diamond,
mark options={solid}
]
table[row sep=crcr]{
8 0.017432381\\
40 5.72e-05\\
80 3.72e-06\\
160 2.35e-07\\
};
\addlegendentry{p=3};

\addplot [
color=black,
solid,
mark=square,
mark options={solid}
]
table[row sep=crcr]{
8 0.003132092\\
40 1.79e-05\\
80 1.17e-06\\
160 7.42e-08\\
};
\addlegendentry{p=4};

\addplot [
color=black,
solid,
mark=o,
mark options={solid}
]
table[row sep=crcr]{
8 0.001439654\\
40 3.82e-07\\
80 7.99e-09\\
160 2.28e-10\\
};
\addlegendentry{p=5};

\end{axis}
\end{tikzpicture}%
}
\subfigure[Relative error in displacement, $u_y$]{
%
%
%
%
\begin{tikzpicture}

\begin{axis}[%
width=\figurewidth,
height=\figureheight,
scale only axis,
xmode=log,
xmin=1,
xmax=1000,
xminorticks=true,
xlabel={log(Number of elements)},
ymode=log,
ymin=1e-10,
ymax=1,
yminorticks=true,
ylabel={$\text{log(Rel. error in displacement (u}_\text{y}\text{))}$},
legend style={draw=black,fill=white,legend cell align=left}
]
\addplot [
color=black,
solid,
mark=triangle,
mark options={solid,,rotate=90}
]
table[row sep=crcr]{
8 0.171196843\\
40 0.006353531\\
80 0.001598334\\
160 0.000400805\\
};
\addlegendentry{p=2};

\addplot [
color=black,
solid,
mark=diamond,
mark options={solid}
]
table[row sep=crcr]{
8 0.019051983\\
40 3.45e-05\\
80 2.18e-06\\
160 1.37e-07\\
};
\addlegendentry{p=3};

\addplot [
color=black,
solid,
mark=square,
mark options={solid}
]
table[row sep=crcr]{
8 0.001936228\\
40 4.67e-06\\
80 3.06e-07\\
160 1.94e-08\\
};
\addlegendentry{p=4};

\addplot [
color=black,
solid,
mark=o,
mark options={solid}
]
table[row sep=crcr]{
8 0.000402525\\
40 1.22e-07\\
80 3.52e-09\\
160 1.35e-10\\
};
\addlegendentry{p=5};

\end{axis}
\end{tikzpicture}%
}
\caption{Crack in an isotropic medium: Convergence in the relative error in the displacement.}
\label{fig:isocrkres}
\end{figure}

In Table \ref{table:isocrksingularity}, the convergence of the order of singularity with increasing number of elements is shown. With decreasing element size, it is seen that the numerically estimated order of singularity converges to the theoretical value 0.5. Also, the influence of the order of the shape functions on the convergence of the order of singularity is shown. As expected, the theoretical value of 0.5 can be achieved by either increasing the order of the shape functions or by decreasing the element size. \fref{fig:isocrkres} shows the relative error in the displacement norm for each component (i.e., $u_x$ and $u_y$) with the increasing number of elements along the boundary. It can be seen that increasing the number of elements leads to decreasing error.

\begin{figure}
\centering
\setlength\figureheight{8cm} 
\setlength\figurewidth{10cm}
\input{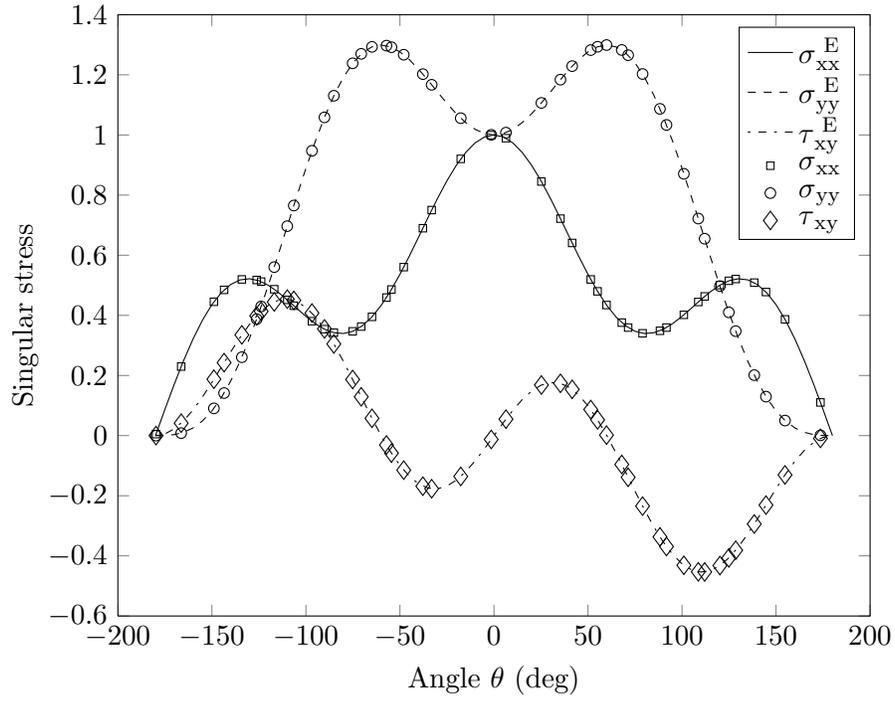}
\caption{Angular distribution of singular stress ahead of a crack tip for a crack in an isotropic medium under mode I loading conditions, i.e., $K_I=$ 1 and $K_{II} = $ 0.}
\label{fig:isomatsingstressD}
\end{figure}

\Eref{eq-stress-bd-slu} is the parametric equation for the stress field in the polar coordinates $r$ and $\theta$. The terms $\left(r_{\eta}^{\lambda_{i}+1}(\eta)\boldsymbol{\psi}_{\sigma i}(\eta)\right)$ in \Eref{eq-stress-bd-slu} together with $\theta(\eta)$ in \Eref{eq-polar-theta} are the stress modes describing the angular stress distribution at a constant radial coordinate $r$. \fref{fig:isomatsingstressD} shows the stress distribution ahead of the crack tip for mode I fracture problem. Each of the stress mode is normalized with its value of $\sigma_{yy}$ at $\theta=$ 0$^\circ$. A total of 8 third order elements were used to discretize the boundary. The results from the scaled boundary formulation are compared with the analytical solutions (represented by the superscript $E$). A very good agreement is observed.

\begin{table}[htpb]
\centering
\caption{Convergence of the order of singularity for an isotropic plate with a crack.}
\begin{tabular}{lrrrr}
\hline 
Number of & \multicolumn{4}{c}{Order of shape functions}\\
\cline{2-5}
Elements & $p=$2 & $p=$ 3 & $p=$ 4 & $p=$ 5 \\
\hline
8 & 0.50826515 & 0.50007949 & 0.49995729 & 0.50000291\\
40 & 0.50035216	 & 0.49999979	 & 0.50000000	& 0.50000000\\
80 & 0.50008862	 & 0.49999999	 & - & - \\	
160 & 0.50002219 & 0.50000000 & - & - \\
\hline
\end{tabular}
\label{table:isocrksingularity}
\end{table}

\begin{figure}[htpb]
\centering
\subfigure[]{\scalebox{0.5}{\input{./Figures/bimatcrk.pstex_t}}}
\subfigure[]{\includegraphics[scale=0.5]{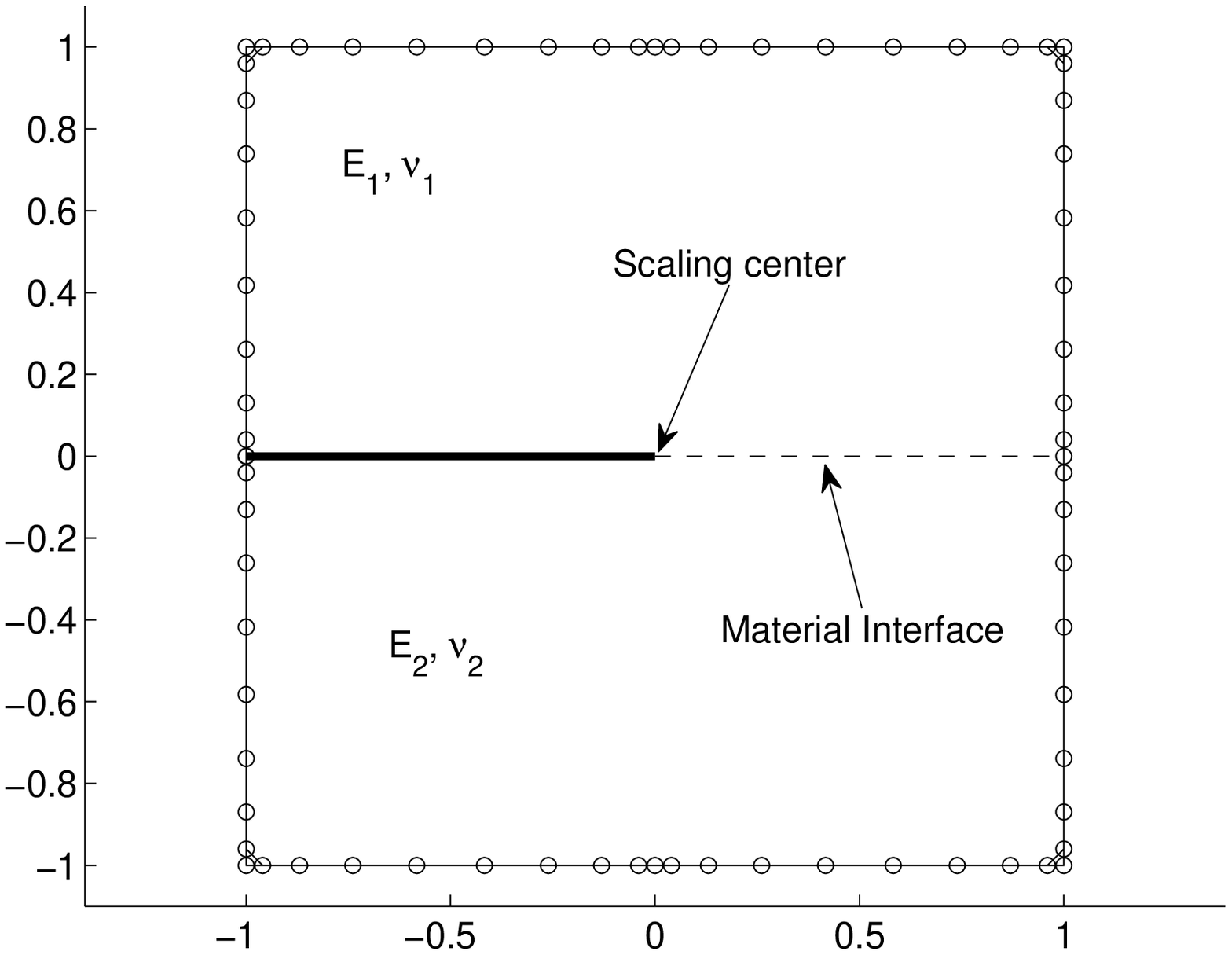}}
\caption{Crack aligned to a bi-material interface: (a) geometry and (b) SBFEM Mesh}
\label{fig:bimatprobdescription}
\end{figure}

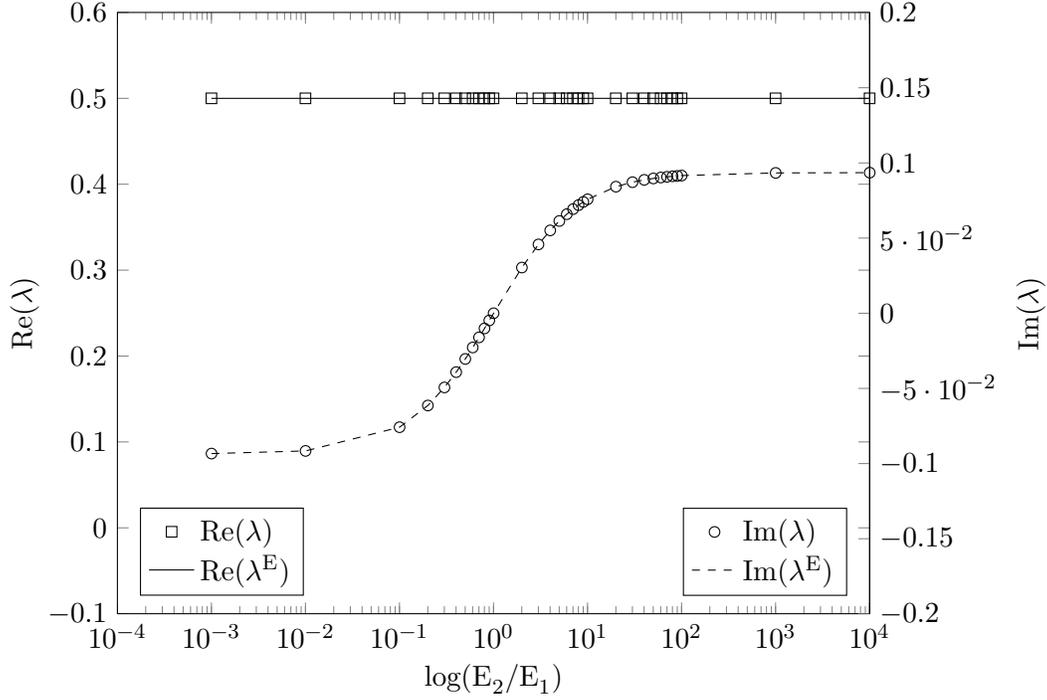
\begin{figure}
\centering
\setlength\figureheight{8cm} 
\setlength\figurewidth{10cm}
%
%
%
%
\begin{tikzpicture}

\begin{axis}[%
width=\figurewidth,
height=\figureheight,
scale only axis,
xmode=log,
xmin=0.0001,
xmax=10000,
xminorticks=true,
xlabel={$\text{log(E}_\text{2}\text{/E}_\text{1}\text{)}$},
ymin=-0.1,
ymax=0.6,
ylabel={$\text{Re(}\lambda\text{)}$},
legend style={draw=black,fill=white,legend cell align=left},
legend pos = south west
]
\addplot [
color=black,
only marks,
mark=square,
mark options={solid}
]
table[row sep=crcr]{
0.001 0.49999978966\\
0.01 0.499999796761\\
0.1 0.499999853267\\
0.2 0.499999894317\\
0.3 0.49999992129\\
0.4 0.499999939294\\
0.5 0.49999995134\\
0.6 0.499999959299\\
0.7 0.499999964387\\
0.8 0.49999996742\\
0.9 0.49999996896\\
1 0.499999967431\\
2 0.49999995134\\
3 0.499999928109\\
4 0.499999909178\\
5 0.499999894317\\
6 0.49999988254\\
7 0.499999873045\\
8 0.499999865254\\
9 0.499999858759\\
10 0.499999853267\\
20 0.499999824827\\
30 0.499999813782\\
40 0.499999807922\\
50 0.499999804291\\
60 0.499999801822\\
70 0.499999800031\\
80 0.499999798678\\
90 0.499999797615\\
100 0.499999853267\\
1000 0.49999979676\\
10000 0.49999978892\\
};
\addlegendentry{$\text{Re(}\lambda\text{)}$};

\addplot [
color=black,
solid
]
table[row sep=crcr]{
0.001 0.5\\
0.01 0.5\\
0.1 0.5\\
0.2 0.5\\
0.3 0.5\\
0.4 0.5\\
0.5 0.5\\
0.6 0.5\\
0.7 0.5\\
0.8 0.5\\
0.9 0.5\\
1 0.5\\
2 0.5\\
3 0.5\\
4 0.5\\
5 0.5\\
6 0.5\\
7 0.5\\
8 0.5\\
9 0.5\\
10 0.5\\
20 0.5\\
30 0.5\\
40 0.5\\
50 0.5\\
60 0.5\\
70 0.5\\
80 0.5\\
90 0.5\\
100 0.5\\
1000 0.5\\
10000 0.5\\
};
\addlegendentry{$\text{Re(}\lambda^{\rm E}\text{)}$};

\end{axis}

\begin{axis}[%
width=\figurewidth,
height=\figureheight,
scale only axis,
xmin=0.0001,
xmax=10000,
xminorticks=true,
ymin=-0.2,
ymax=0.2,
xmode=log,
hide x axis,
axis y line*=right,
ylabel={$\text{Im(}\lambda\text{)}$},
ylabel near ticks,
legend style={draw=black,fill=white,legend cell align=left},
legend pos = south east
]

\addplot [
color=black,
only marks,
mark=o,
mark options={solid}
]
table[row sep=crcr]{
0.001 -0.093350060333\\
0.01 -0.091590349707\\
0.1 -0.075810772121\\
0.2 -0.061379288225\\
0.3 -0.049362051545\\
0.4 -0.039172791529\\
0.5 -0.030406996603\\
0.6 -0.022774916373\\
0.7 -0.01606265453\\
0.8 -0.01010834406\\
0.9 -0.004786914902\\
1 0\\
2 0.030406996603\\
3 0.04578543325\\
4 0.055111574749\\
5 0.061379288225\\
6 0.065885429925\\
7 0.069282110485\\
8 0.071934716638\\
9 0.074063881077\\
10 0.075810772122\\
20 0.084193100862\\
30 0.087194551611\\
40 0.088737727537\\
50 0.089677723112\\
60 0.090310379969\\
70 0.090765252631\\
80 0.091108049405\\
90 0.091375649146\\
100 0.09159034971\\
1000 0.093350060329\\
10000 0.093528080505\\
};
\addlegendentry{$\text{Im(}\lambda\text{)}$};

\addplot [
color=black,
dashed
]
table[row sep=crcr]{
0.001 -0.093351326688\\
0.01 -0.091591589141\\
0.1 -0.07581177769\\
0.2 -0.061380090186\\
0.3 -0.049362689964\\
0.4 -0.039173294675\\
0.5 -0.030407385334\\
0.6 -0.022775206626\\
0.7 -0.01606285883\\
0.8 -0.010108472473\\
0.9 -0.004786975677\\
1 0\\
2 0.030407385334\\
3 0.04578602387\\
4 0.055110858695\\
5 0.061380090186\\
6 0.065886294536\\
7 0.069283022858\\
8 0.071935666648\\
9 0.074064861521\\
10 0.075811778\\
20 0.084194229039\\
30 0.087195724567\\
40 0.088738923708\\
50 0.089678933487\\
60 0.090311599933\\
70 0.090766479504\\
80 0.091109281491\\
90 0.091376885308\\
100 0.091591589141\\
1000 0.093351326688\\
10000 0.093529349587\\
};
\addlegendentry{$\text{Im(}\lambda^{\rm E}\text{)}$};

\end{axis}
\end{tikzpicture}%
\caption{Bimaterial interface crack: order of singularity $(\lambda)$ as a function of the Young's modulus mismatch.}
\label{fig:bimatsing}
\end{figure}

\begin{figure}
\centering
\setlength\figureheight{8cm} 
\setlength\figurewidth{10cm}
\input{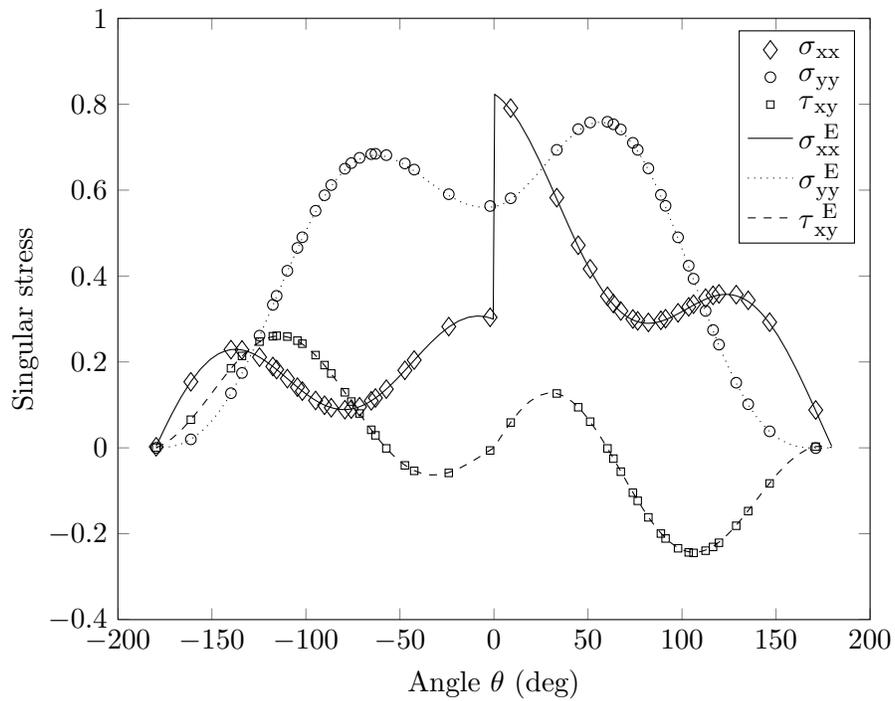}
\caption{Angular distribution of singular stress ahead of a crack tip for a crack aligned to the bimaterial interface. The ratio of Young's modulus is $E_1/E_2=$ 10. The singular stress fields are obtained for a crack length $a$ = 0.5 and $K_I/p\sqrt{\pi a}=$ 1.123 and $K_{II}/p\sqrt{\pi a}=$ -0.123.}
\label{fig:bimatsingstressD}
\end{figure}

\paragraph{Crack aligned to a bi-material interface} The second test case is a crack aligned to the bimaterial interface. \fref{fig:bimatprobdescription} shows the geometry and the corresponding finite element mesh. Again, only the boundary of the domain is discretized. The Lagrange shape functions of order $p=$ 9 and a total of eight elements are used to discretize the domain in the circumferential direction. The order of the singularity computed from the semi-analytical approach is compared with the analytical expression  given by, $0.5 \pm i \epsilon$, where $\epsilon$ is the oscillatory index given by:
\begin{equation}
\epsilon = \frac{1}{2\pi} \mathrm{log} \left( \frac{1 - \beta}{1+\beta} \right)
\label{eqn:osciindex}
\end{equation}
where $\beta$ is the second Dundur's parameter
\begin{equation}
\beta = \frac{\mu_1(\kappa_2-1) + \mu_2(\kappa_1-1)}{\mu_1(\kappa_2+1) + \mu_2(\kappa_1+1)}
\end{equation}
and,
\begin{equation}
\kappa_i = \left\{ \begin{array}{cr} (3-\nu_i)/(1+\nu_i) & \textup{plane stress} \\ 3 -4\nu_i & \textup{plane strain} \end{array} \right.
\end{equation}
where $\mu_i, \nu_i$ and $\kappa_i, (i=1,2)$ are the shear modulus, Poisson's ratio and the Kolosov coefficient, respectively. \fref{fig:bimatsing} shows the influence of the ratio of Young's modulus on the order of the singularity. The angular stress distribution ahead of the crack aligned to a bimaterial interface is shown in \fref{fig:bimatsingstressD}.  The ratio of Young's modulus is $E_1/E_2=$ 10. The singular stress fields are obtained for a crack length $a$ = 0.5 and $K_I/p\sqrt{\pi a}=$ 1.123 and $K_{II}/p\sqrt{\pi a}=$ -0.123. A very good agreement is observed.

\begin{figure}[htpb]
\centering
\subfigure[]{\scalebox{0.5}{\input{./Figures/trimatnocrk.pstex_t}}}
\subfigure[]{\includegraphics[scale=0.5]{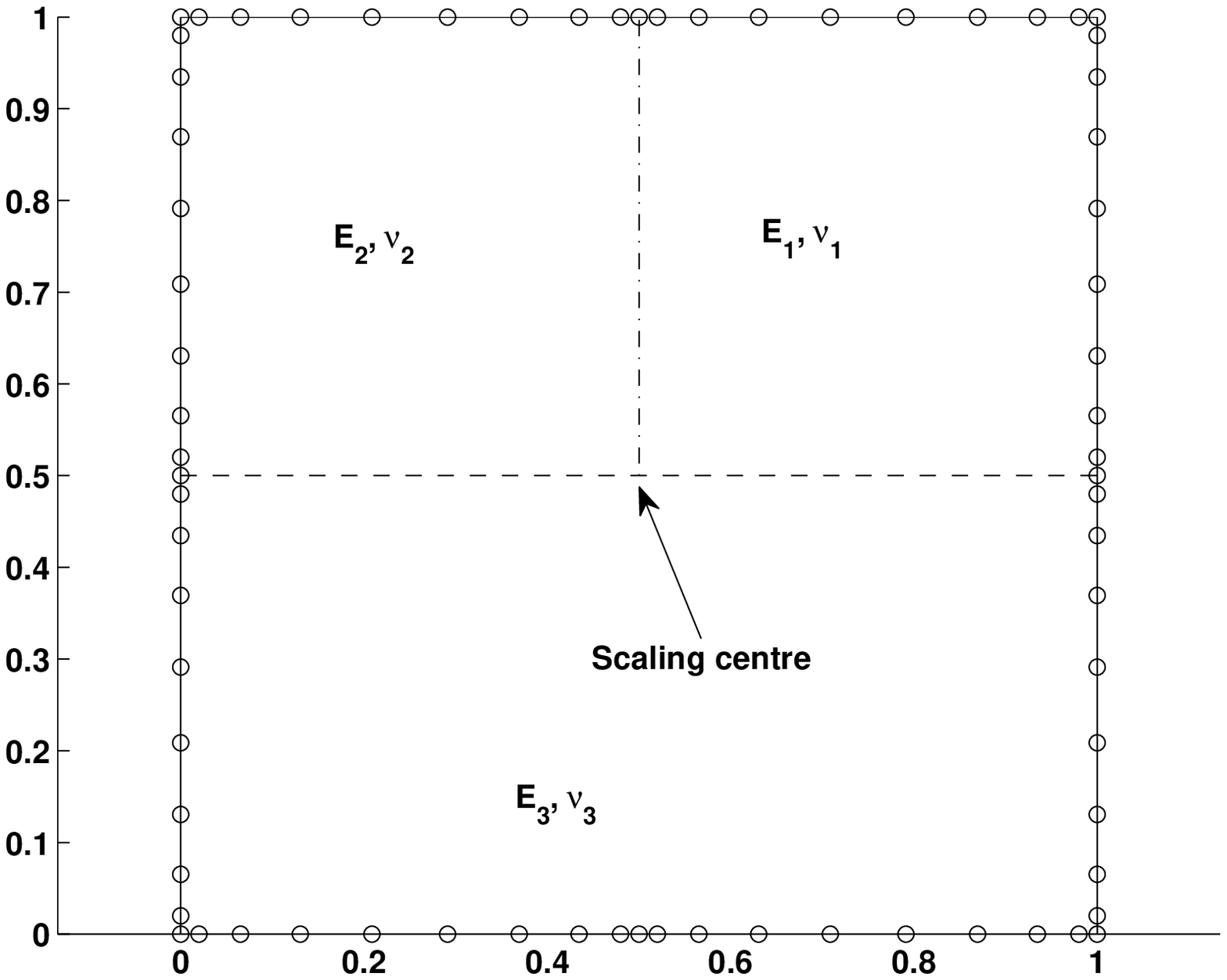}}
\subfigure[]{\scalebox{0.5}{\input{./Figures/trimatcrk.pstex_t}}}
\subfigure[]{\includegraphics[scale=0.5]{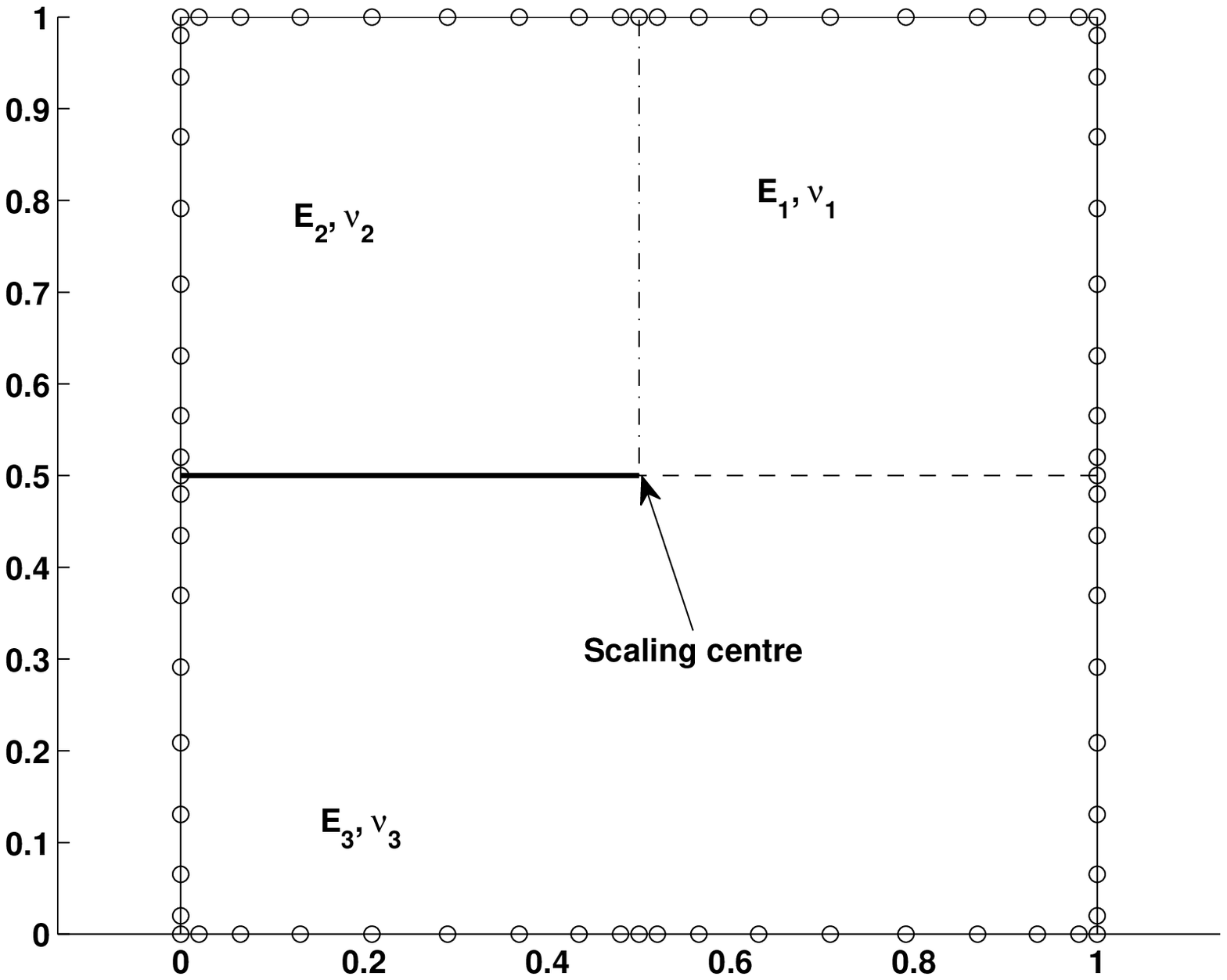}}
\caption{(a-b) a fully bonded triple junction and SBFEM mesh and (c-d) a triple junction with a crack aligned to one of the interfaces.}
\label{fig:crkatint}
\end{figure}

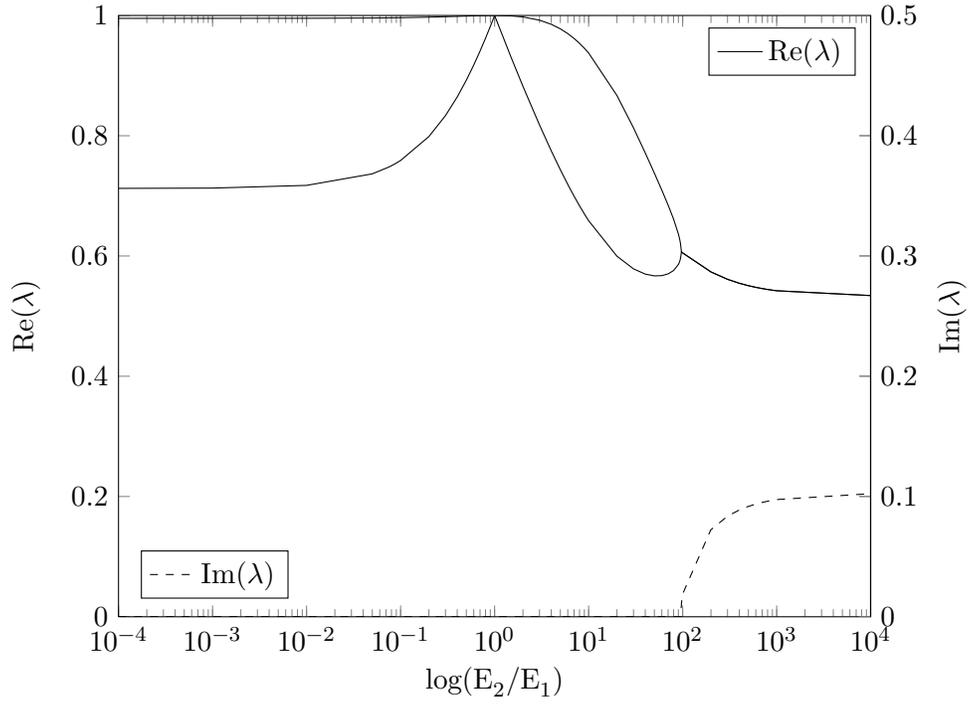
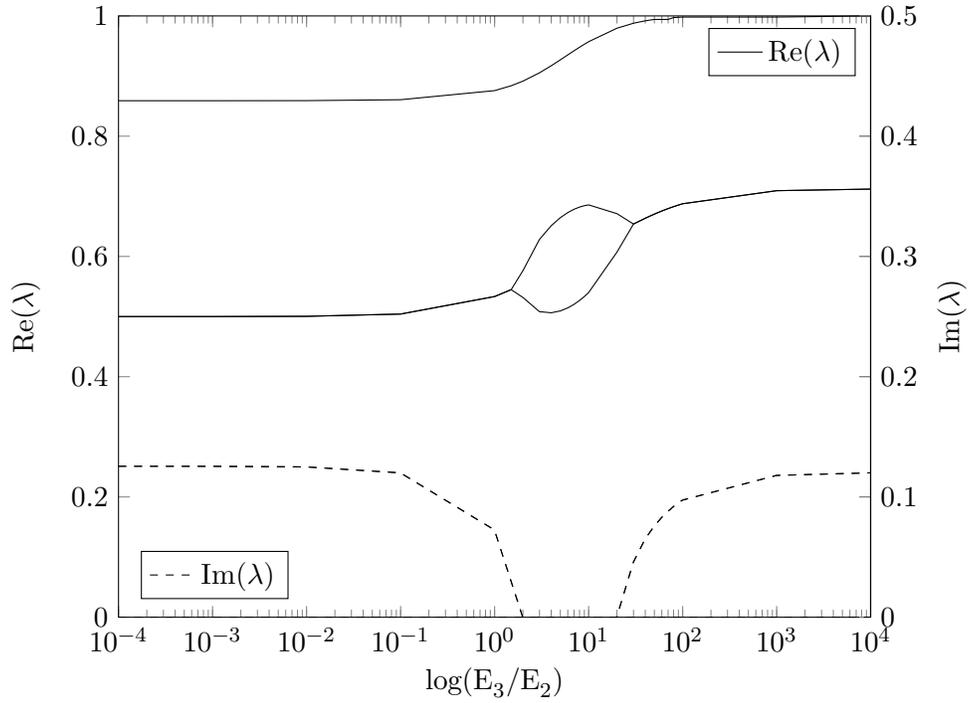
\begin{figure}
\centering
\setlength\figureheight{8cm} 
\setlength\figurewidth{10cm}
\subfigure[$E_3/E_1=$ 10, $E_1=$ 1]{
%
%
%
%
\begin{tikzpicture}

\begin{axis}[%
width=\figurewidth,
height=\figureheight,
scale only axis,
xmode=log,
xmin=0.0001,
xmax=10000,
xminorticks=true,
xlabel={$\text{log(E}_\text{2}\text{/E}_\text{1}\text{)}$},
ymin=0,
ymax=1,
ylabel={$\text{Re(}\lambda\text{)}$},
legend style={draw=black,fill=white,legend cell align=left}
]
\addplot [
color=black,
solid
]
table[row sep=crcr]{
0.0001 0.995225636\\
0.001 0.995236291\\
0.01 0.995341902\\
0.05 0.995791261\\
0.075 0.996056094\\
0.08 0.996107632\\
0.085 0.9961587\\
0.09 0.996209304\\
0.095 0.996259445\\
0.1 0.996309128\\
0.2 0.99721118\\
0.3 0.997953978\\
0.4 0.998556943\\
0.5 0.999036372\\
0.6 0.999406068\\
0.7 0.99967781\\
0.8 0.999861726\\
0.9 0.99996658\\
1 1\\
1.5 0.999303543\\
2 0.997555155\\
3 0.992163544\\
4 0.985342869\\
5 0.977783897\\
6 0.969849638\\
7 0.961749454\\
8 0.953611838\\
9 0.945518932\\
10 0.937524519\\
20 0.866885126\\
30 0.813076769\\
40 0.77139761\\
50 0.737814311\\
60 0.709578797\\
70 0.684706902\\
80 0.661402214\\
90 0.636683281\\
91 0.633820686\\
92 0.63078399\\
93 0.627498817\\
94 0.623828101\\
95 0.619456815\\
96 0.613110574\\
97 0.606850592\\
98 0.60626004\\
99 0.605678965\\
100 0.605107138\\
200 0.573471011\\
300 0.561203588\\
400 0.554676782\\
500 0.550623472\\
600 0.547861126\\
700 0.545857568\\
800 0.544337847\\
900 0.543145554\\
1000 0.542185148\\
10000 0.534182239\\
};
\addlegendentry{$\text{Re(}\lambda\text{)}$};

\addplot [
color=black,
solid
]
table[row sep=crcr]{
0.0001 0.712457279\\
0.001 0.71291217\\
0.01 0.717419755\\
0.05 0.736600644\\
0.075 0.747943448\\
0.08 0.75015725\\
0.085 0.752353516\\
0.09 0.75453254\\
0.095 0.756694602\\
0.1 0.758839981\\
0.2 0.798612567\\
0.3 0.833444003\\
0.4 0.864429288\\
0.5 0.892314251\\
0.6 0.91763339\\
0.7 0.940785249\\
0.8 0.962076575\\
0.9 0.981749654\\
1 1\\
1.5 0.930048434\\
2 0.882282777\\
3 0.817996152\\
4 0.774957764\\
5 0.743386875\\
6 0.718977854\\
7 0.699440872\\
8 0.683410495\\
9 0.6700077\\
10 0.658634642\\
20 0.599644185\\
30 0.578598482\\
40 0.569884118\\
50 0.566906479\\
60 0.567309702\\
70 0.570247915\\
80 0.575804503\\
90 0.585857986\\
91 0.587392676\\
92 0.589124076\\
93 0.591125983\\
94 0.593534908\\
95 0.596665339\\
96 0.601791139\\
97 0.606850592\\
98 0.60626004\\
99 0.605678965\\
100 0.605107138\\
200 0.573471011\\
300 0.561203588\\
400 0.554676782\\
500 0.550623472\\
600 0.547861126\\
700 0.545857568\\
800 0.544337847\\
900 0.543145554\\
1000 0.542185148\\
10000 0.534182239\\
};

\end{axis}

\begin{axis}[%
width=\figurewidth,
height=\figureheight,
scale only axis,
xmin=0.0001,
xmax=10000,
xminorticks=true,
ymin=0,
ymax=0.5,
xmode=log,
hide x axis,
axis y line*=right,
ylabel={$\text{Im(}\lambda\text{)}$},
ylabel near ticks,
legend style={draw=black,fill=white,legend cell align=left},
legend pos = south west
]
\addplot [
color=black,
dashed
]
table[row sep=crcr]{
0.0001 0\\
0.001 0\\
0.01 0\\
0.05 0\\
0.1 0\\
1 0\\
2 0\\
3 0\\
4 0\\
5 0\\
6 0\\
7 0\\
8 0\\
9 0\\
10 0\\
20 0\\
30 0\\
40 0\\
50 0\\
60 0\\
70 0\\
80 0\\
90 0\\
91 0\\
92 0\\
93 0\\
94 0\\
95 0\\
96 0\\
97 0.008007142\\
98 0.012594381\\
99 0.015859373\\
100 0.018516394\\
200 0.072180656\\
300 0.083446977\\
400 0.088641589\\
500 0.091638799\\
600 0.093590122\\
700 0.094961688\\
800 0.095978412\\
900 0.096762191\\
1000 0.097384831\\
10000 0.10228642\\
};
\addlegendentry{$\text{Im(}\lambda\text{)}$};

\end{axis}

\end{tikzpicture}%
}
\subfigure[$E_1/E_2=$ 10, $E_2=$ 1]{
%
%
%
%
\begin{tikzpicture}

\begin{axis}[%
width=\figurewidth,
height=\figureheight,
scale only axis,
xmode=log,
xmin=0.0001,
xmax=10000,
xminorticks=true,
xlabel={$\text{log(E}_\text{3}\text{/E}_\text{2}\text{)}$},
ymin=0,
ymax=1,
ylabel={$\text{Re(}\lambda\text{)}$},
legend style={draw=black,fill=white,legend cell align=left}
]
\addplot [
color=black,
solid
]
table[row sep=crcr]{
0.0001 0.500004313\\
0.001 0.500043129\\
0.01 0.500430172\\
0.1 0.504192784\\
1 0.533325052\\
1.5 0.544789582\\
2 0.576353416\\
3 0.628353492\\
4 0.651113374\\
5 0.664492906\\
6 0.672971872\\
7 0.6785091\\
8 0.682129817\\
9 0.684424463\\
10 0.685754793\\
20 0.671163074\\
30 0.653867689\\
40 0.66352676\\
50 0.670402204\\
60 0.675560429\\
70 0.679578981\\
80 0.682800508\\
90 0.685441987\\
100 0.687647817\\
1000 0.709375548\\
10000 0.712028461\\
};
\addlegendentry{$\text{Re(}\lambda\text{)}$};

\addplot [
color=black,
solid
]
table[row sep=crcr]{
0.0001 0.500004313\\
0.001 0.500043129\\
0.01 0.500430172\\
0.1 0.504192784\\
1 0.533325052\\
1.5 0.544789582\\
2 0.531789685\\
3 0.508141789\\
4 0.506329765\\
5 0.509404928\\
6 0.514453192\\
7 0.520411038\\
8 0.526797078\\
9 0.533370253\\
10 0.540004796\\
20 0.607142738\\
30 0.653867689\\
40 0.66352676\\
50 0.670402204\\
60 0.675560429\\
70 0.679578981\\
80 0.682800508\\
90 0.685441987\\
100 0.687647817\\
1000 0.709375548\\
10000 0.712028461\\
};

\addplot [
color=black,
solid
]
table[row sep=crcr]{
0.0001 0.858962126\\
0.001 0.858977313\\
0.01 0.859129206\\
0.1 0.860650118\\
1 0.875792043\\
1.5 0.883880393\\
2 0.891555102\\
3 0.905391387\\
4 0.917132393\\
5 0.92696516\\
6 0.93518614\\
7 0.942087619\\
8 0.947920126\\
9 0.952886937\\
10 0.957149484\\
20 0.97956082\\
30 0.987851627\\
40 0.991894696\\
50 0.994188603\\
60 0.994188603\\
70 0.994188603\\
80 0.997252174\\
90 0.99774324\\
100 0.998112922\\
1000 0.998112922\\
10000 0.999999885\\
};

\end{axis}

\begin{axis}[%
width=\figurewidth,
height=\figureheight,
scale only axis,
xmin=0.0001,
xmax=10000,
xminorticks=true,
ymin=0,
ymax=0.5,
xmode=log,
hide x axis,
axis y line*=right,
ylabel={$\text{Im(}\lambda\text{)}$},
ylabel near ticks,
legend style={draw=black,fill=white,legend cell align=left},
legend pos = south west
]
\addplot [
color=black,
dashed
]
table[row sep=crcr]{
0.0001 0.125481361\\
0.001 0.125431617\\
0.01 0.124934541\\
0.1 0.119997589\\
1 0.072552358\\
2 0\\
3 0\\
4 0\\
5 0\\
6 0\\
7 0\\
8 0\\
9 0\\
10 0\\
20 0\\
30 0.045600324\\
40 0.064847757\\
50 0.075668051\\
60 0.082827617\\
70 0.087962825\\
80 0.091840538\\
90 0.094877843\\
100 0.09732368\\
1000 0.117855961\\
10000 0.119995829\\
};
\addlegendentry{$\text{Im(}\lambda\text{)}$};

\addplot [
color=black,
dashed,
forget plot
]
table[row sep=crcr]{
0.0001 0.125481361\\
0.001 0.125431617\\
0.01 0.124934541\\
0.1 0.119997589\\
1 0.072552358\\
2 0\\
3 0\\
4 0\\
5 0\\
6 0\\
7 0\\
8 0\\
9 0\\
10 0\\
20 0\\
30 0.045600324\\
40 0.064847757\\
50 0.075668051\\
60 0.082827617\\
70 0.087962825\\
80 0.091840538\\
90 0.094877843\\
100 0.09732368\\
1000 0.117855961\\
10000 0.119995829\\
};
\addplot [
color=black,
dashed,
forget plot
]
table[row sep=crcr]{
0.0001 0\\
0.001 0\\
0.01 0\\
0.1 0\\
1 0\\
2 0\\
3 0\\
4 0\\
5 0\\
6 0\\
7 0\\
8 0\\
9 0\\
10 0\\
20 0\\
30 0\\
40 0\\
50 0\\
60 0\\
70 0\\
80 0\\
90 0\\
100 0\\
1000 0\\
10000 1.34764e-07\\
};
\end{axis}

\end{tikzpicture}%
}
\caption{Order of singularities for a triple junction: (a) a fully bonded junction and (b) a crack along the interface between the material 2 and material 3.}
\label{fig:tjunsingular}
\end{figure}

\paragraph{Triple junction} As a last example, we consider the two specific cases of a triple junction: (a) a fully bonded triple junction and (b) a triple junction with a crack aligned to one of the interfaces. \fref{fig:crkatint} shows the problem description and the corresponding SBFEM mesh. Note that only the boundary of the domain needs to be discretized with 8 elements (on the boundary) with Lagrange shape functions of order $p=$ 9 are used. All materials are assumed to be homogeneous and isotropic. The computed order of singularity for various combination of material parameters are shown in \fref{fig:tjunsingular}. The results shown in \fref{fig:tjunsingular} indicate the existence of significant singularities over the full range of $E_2/E_1$. The singularities are more severe for higher values of $E_2/E_1$ than for lower values. \fref{fig:tjunsingular} also shows the case when the crack is present at the interface of material 3 and material 2. Again, a full range of the order of singularity is shown for various combinations of material parameters. For certain values of $E_3/E_2$, multiple real roots are obtained. The present formulation is efficient in capturing the different orders of singularity. 

\subsection{Crack aligned to the interface}

\subsubsection{Edge Crack aligned to a bimaterial interface in tension} 
\begin{figure}[htpb]
\centering
\scalebox{0.7}{\input{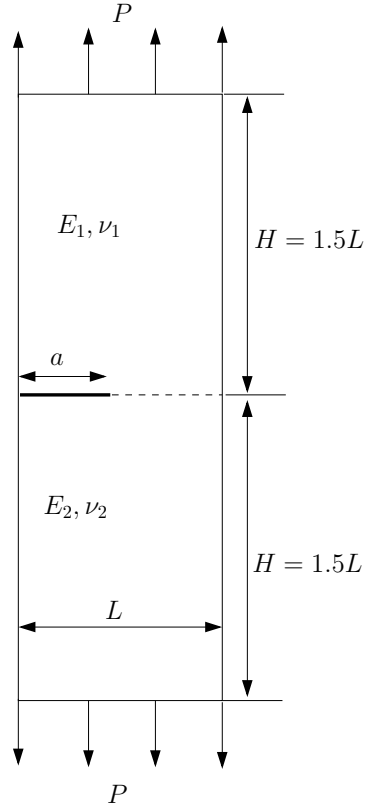}}
\caption{Edge crack aligned to a bimaterial inteface in tension: geometry and boundary conditions.}
\label{fig:crkbimat}
\end{figure}

\begin{figure}[htpb]
\centering
\subfigure[Discretization for a bimaterial interface crack]{\includegraphics[scale=0.5]{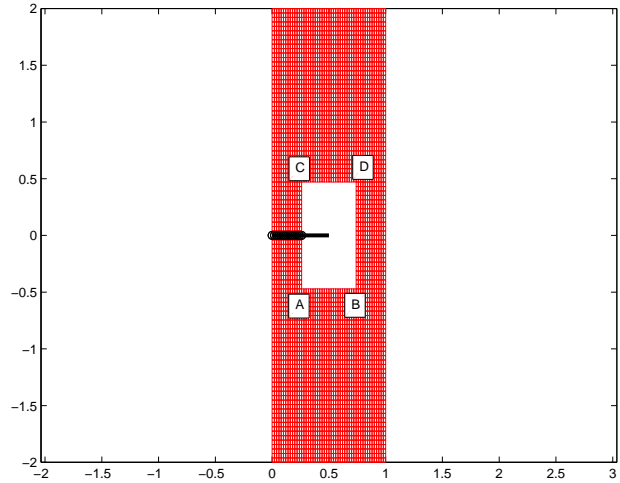}}
\subfigure[Zoomed in view of the region ABCD]{\includegraphics[scale=0.5]{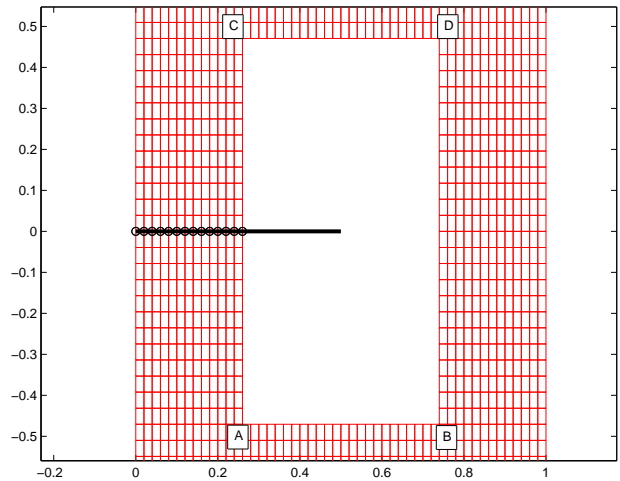}}
\caption{A typical discretization for a crack aligned to a bimaterial interface crack. The region ABCD is identified  from the FE discretization and then replaced with the SBFEM domain. Note that, the SBFEM domain is constructed from the underlying FE discretization. The dofs of nodes within the SBFEM are suppressed during the solution. In this example, the interface aligns with the element edge. To represent the jump across the crack face, the `circled' nodes are enriched with Heaviside function~\cite{sukumarhuang2004}.}
\label{fig:meshdetailsintercrck}
\end{figure}

Consider a bimaterial plate with an edge crack $(a/L=$ 0.5$)$ subjected to uniform far field tension $(P=1)$. The geometry, loading and boundary conditions are shown in \fref{fig:crkbimat}. The material is assumed to be isotropic with Young's modulus $E_1/E_2=$ 2 and Poisson's ratio $\nu_1 = \nu_2 = $ 0.3. A state of plane strain is considered in this example. The discontinuities (both the crack and the material interface) are represented independent of the underlying finite element (FE) mesh. From the underlying FE mesh, the scaled boundary region $\Omega^{\rm sbfem}$ is identified by 
\begin{itemize}
\item Selecting the nodes whose nodal support is intersected by the discontinuous surface, viz., crack or the material interface. The nodes whose support are intersected by the material discontinuity are treated within the XFEM framework.
\item The elements containing these nodes are then selected. Only the information on the boundary of the domain is retained. The dofs of nodes within the domain are condensed during the solution process.
\end{itemize}
Within this region, i.e, $\Omega^{\rm sbfem}$, the scaled boundary formulation is employed to compute the stiffness matrix. Note that no special numerical integration technique is required. The main advantage of this approach is that it does not require a priori knowledge of the asymptotic fields. \fref{fig:meshdetailsintercrck} shows a typical finite element mesh used for this study. The stiffness matrix of the region $A-B-C-D$ is computed by employing the scaled boundary formulation.

Tables \ref{table:mode12meshconve} - \ref{table:mode12layerconve} shows the convergence of the mode I and mode II stress intensity factors with mesh refinement and with number of layers of scaled boundary region in the vicinity of the crack tip. From Table \ref{table:mode12meshconve} it can be seen that with mesh refinement, the numerical SIF approach a constant value. For this study, 3 layers of elements around the crack tip is replaced with the scaled boundary region. The influence of number of layers on the numerical SIFs is demonstrated in Table \ref{table:mode12layerconve} for a structured mesh of 51 $\times$ 102 quadrilateral elements. It is noted that with increasing number of layers around the crack tip, the accuracy of the SIFs improve and approach the values reported in the literature based on the boundary element method~\cite{matsumtotanaka2000} and the XFEM~\cite{liuxiao2004}. It is noted that in~\cite{liuxiao2004}, first 11 terms of the asymptotic solution were retained. From this study, we can conclude that for a given mesh, 4-5 layers of elements should be replaced with the SBFEM domain for reasonable accuracy. However, for more accurate results, more number of layers can be replaced and represented by the SBFEM.

\begin{table}[htpb]
\centering
\caption{Edge Crack aligned to a bimaterial interface in tension: Convergence of the mode I and mode II SIFs with mesh refinement. The number of SBFEM layers around the crack tip is 3.}
\begin{tabular}{llrr}
\hline 
Number & Number & \multicolumn{2}{c}{Stress Intensity Factors}\\
\cline{3-4}
of elements & of dofs & $K_{\rm I}$ & $K_{\rm II}$ \\
\hline
11$\times$22 & 459 & 2.7497 & -0.2700 \\
21$\times$42 & 1769 & 2.7630	 & -0.2729 \\
31$\times$62 & 3879 & 2.7666 & -0.2733 \\
41$\times$82 & 6789 & 2.7679	 & -0.2732 \\
51$\times$102 & 10499 & 2.7685 & -0.2730 \\
\hline
\label{table:mode12meshconve}
\end{tabular}
\end{table}

\begin{table}[htpb]
\centering
\caption{Edge Crack aligned to a bimaterial interface in tension: Convergence of the mode I and mode II SIFs with number of layers of scaled boundary region in the vicinity of the crack tip. A structured quadrilateral mesh (51 $\times$ 102 elements) was used.}
\begin{tabular}{lrr}
\hline 
Number & \multicolumn{2}{c}{Stress Intensity Factors}\\
\cline{2-3}
of layers & $K_{\rm I}$ & $K_{\rm II}$ \\
\hline
3 & 2.7685 & -0.2730 \\
4 & 2.8280 & -0.2707 \\
5 & 2.8219 & -0.2681 \\
6 & 2.8202 & -0.2674 \\
7 & 2.8181 & -0.2671 \\
8 & 2.8170 & -0.2670 \\
Ref.~\cite{matsumtotanaka2000} & 2.8190 & -0.2680 \\
Ref.~\cite{liuxiao2004} & 2.8440 & -0.2670 \\
\hline
\label{table:mode12layerconve}
\end{tabular}
\end{table}

\begin{figure}[htpb]
\centering
\subfigure[]{\includegraphics[scale=0.4]{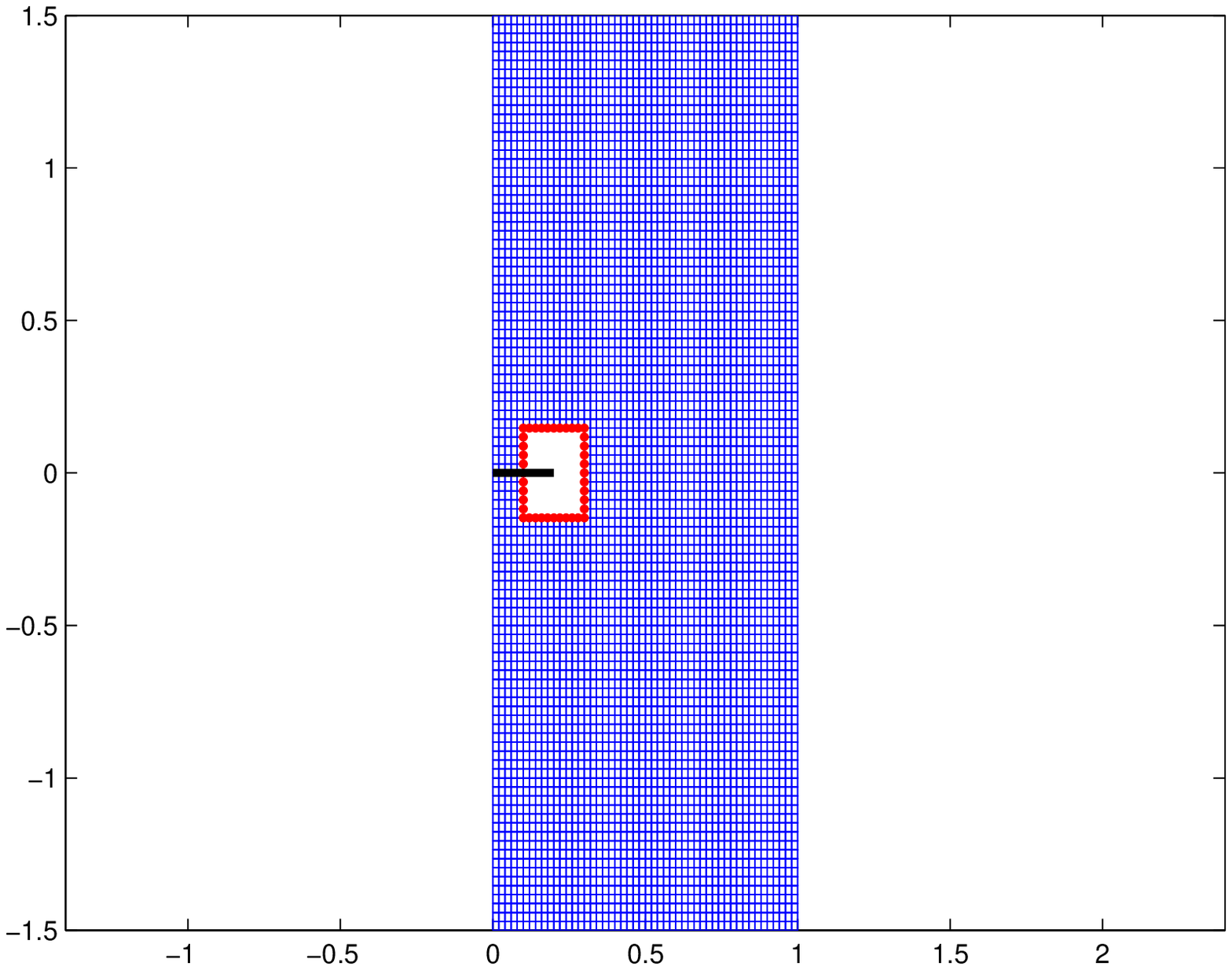}}
\subfigure[]{\includegraphics[scale=0.4]{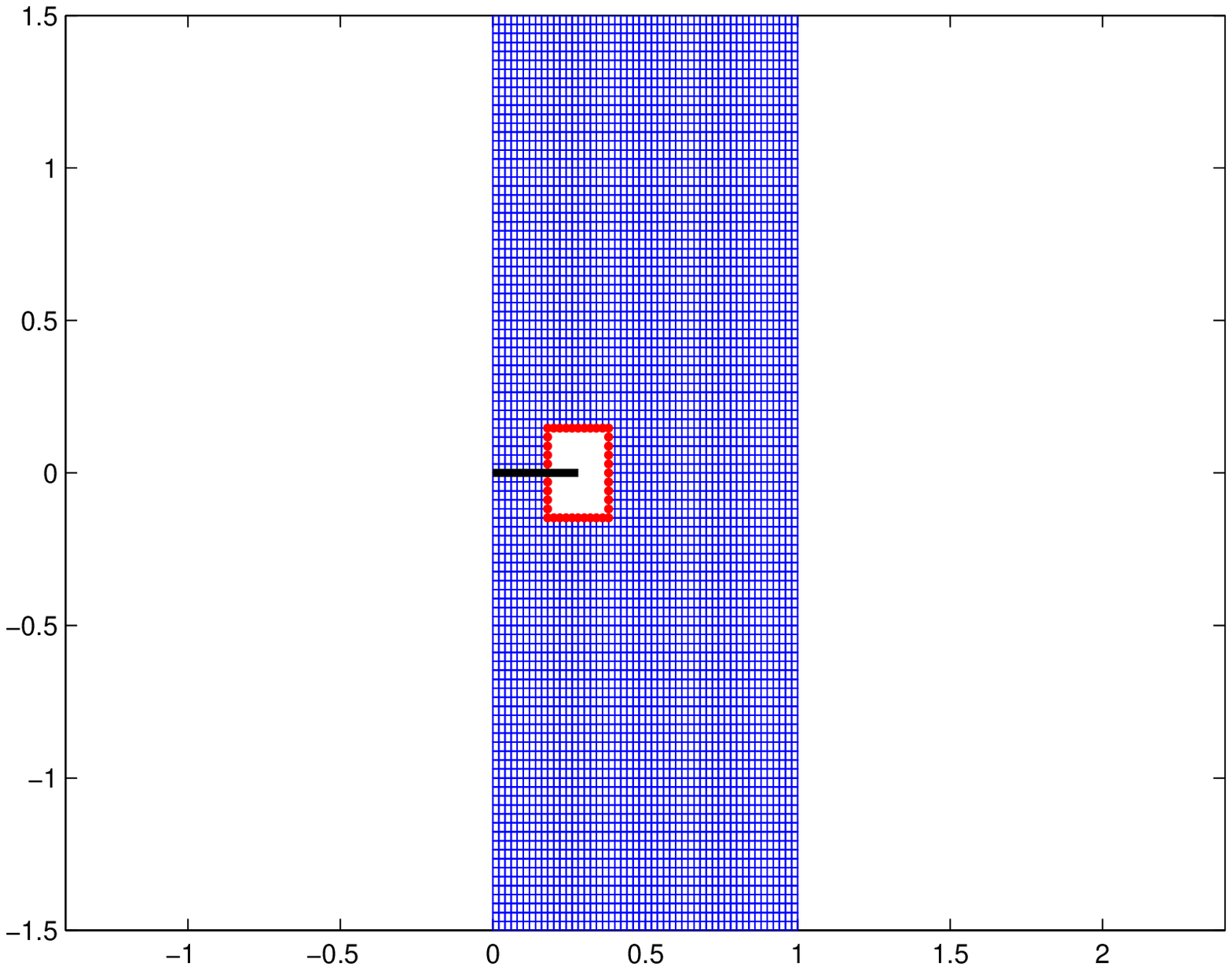}}
\subfigure[]{\includegraphics[scale=0.4]{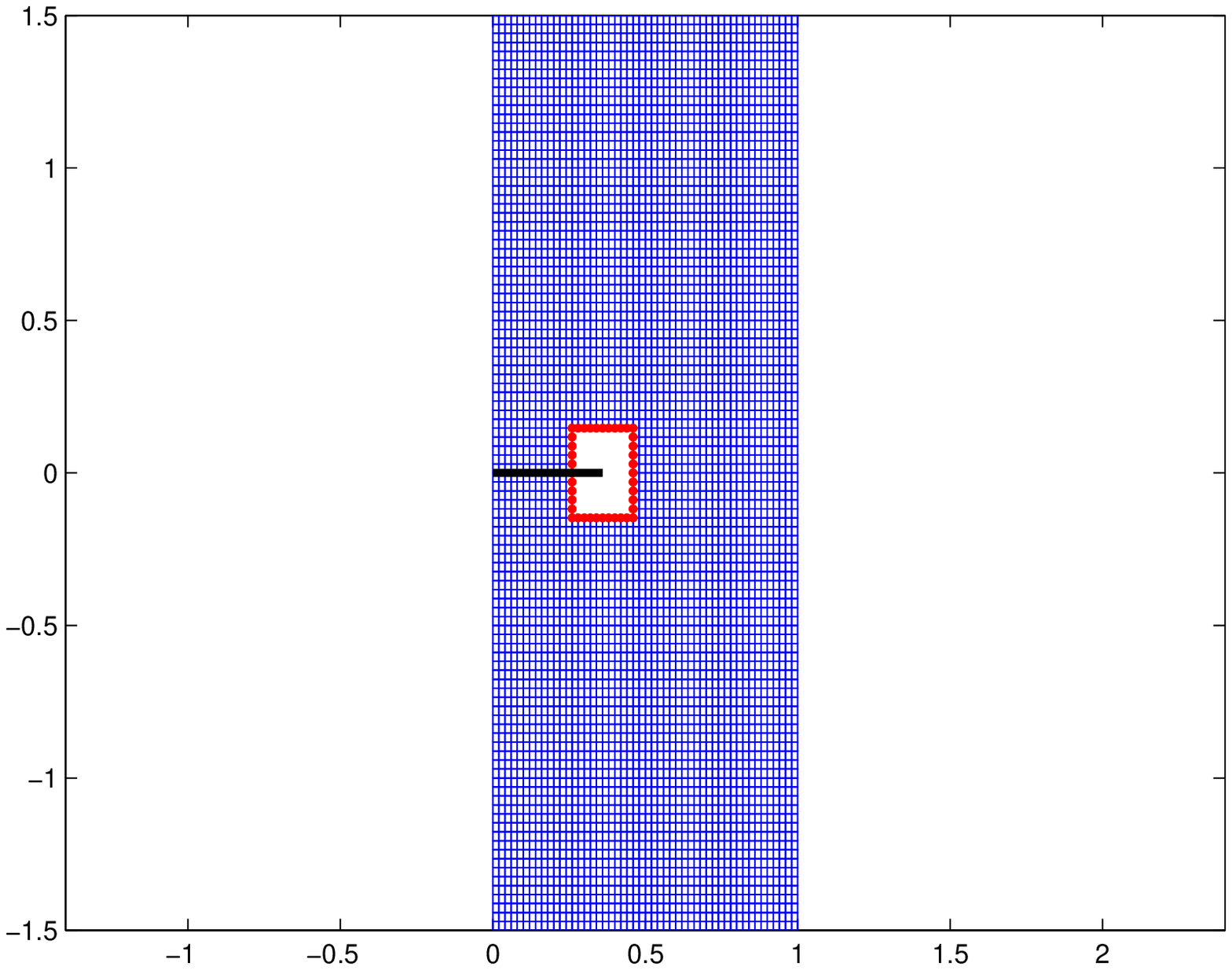}}
\subfigure[]{\includegraphics[scale=0.4]{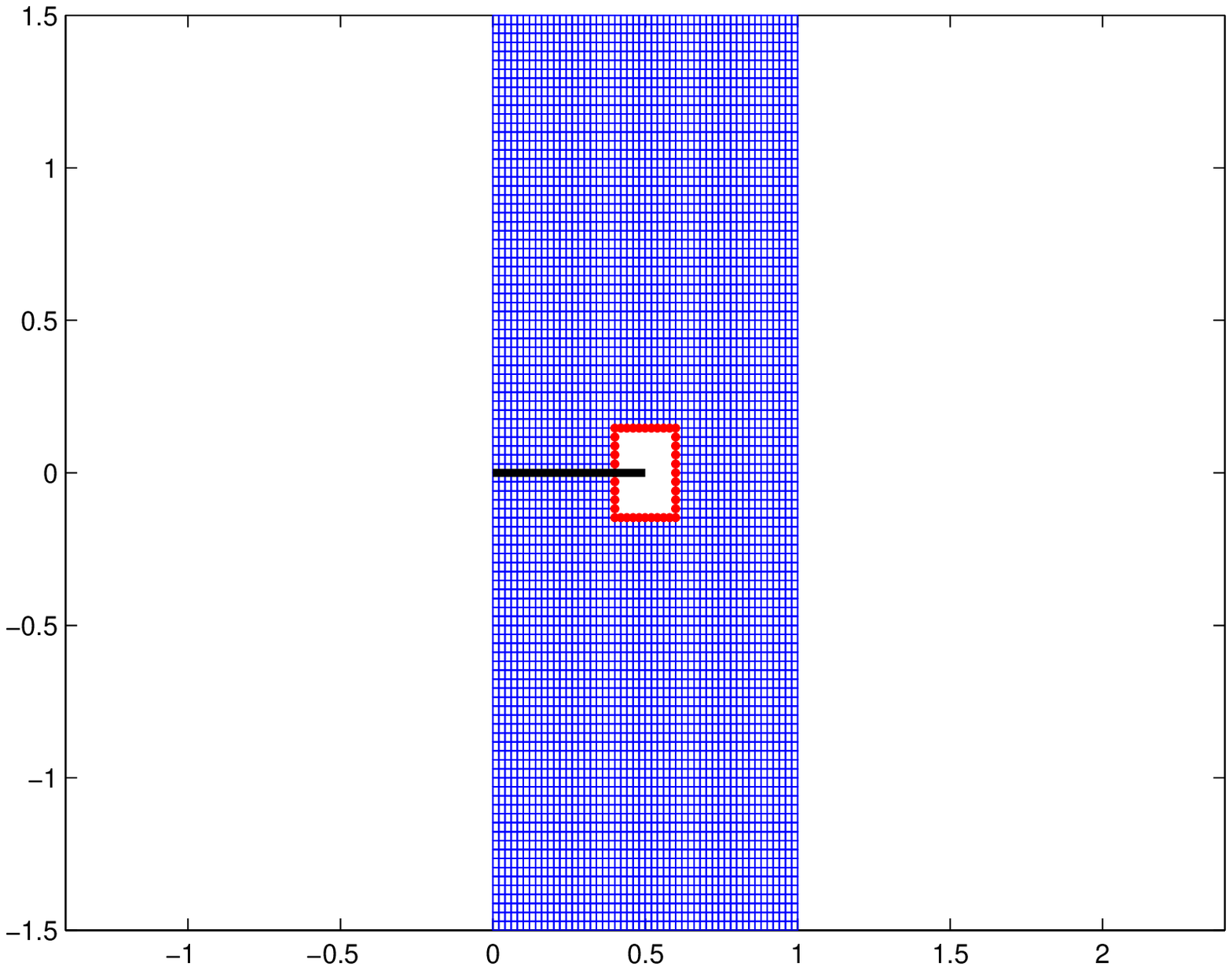}}
\caption{Edge Crack aligned to a bimaterial interface in tension: crack growth along the interface. Meshes at the initial step and at intermediate steps.}
\label{fig:bimatcrkmesh}
\end{figure}

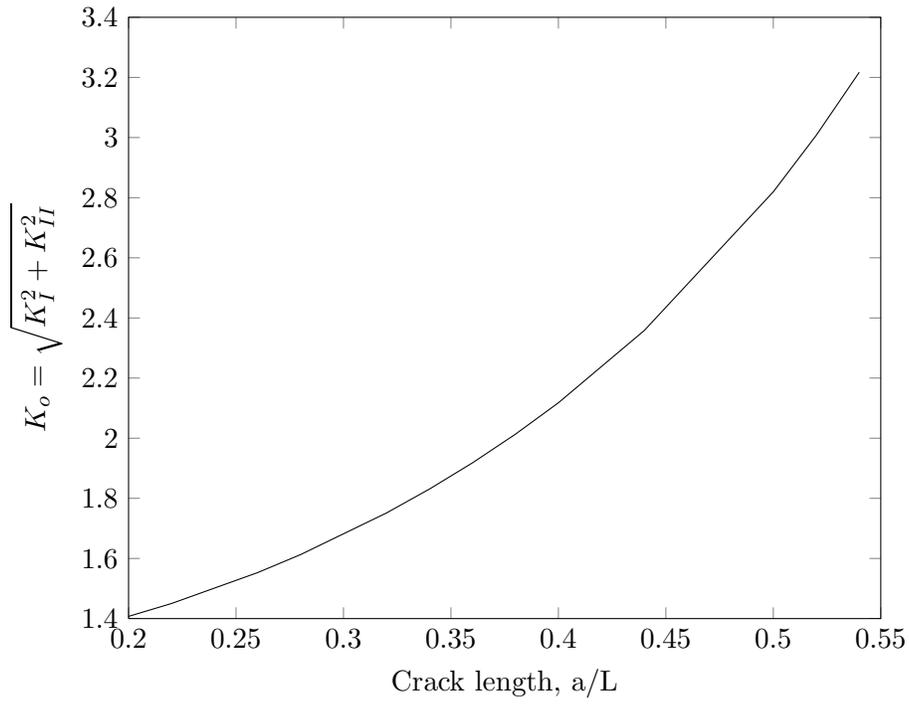
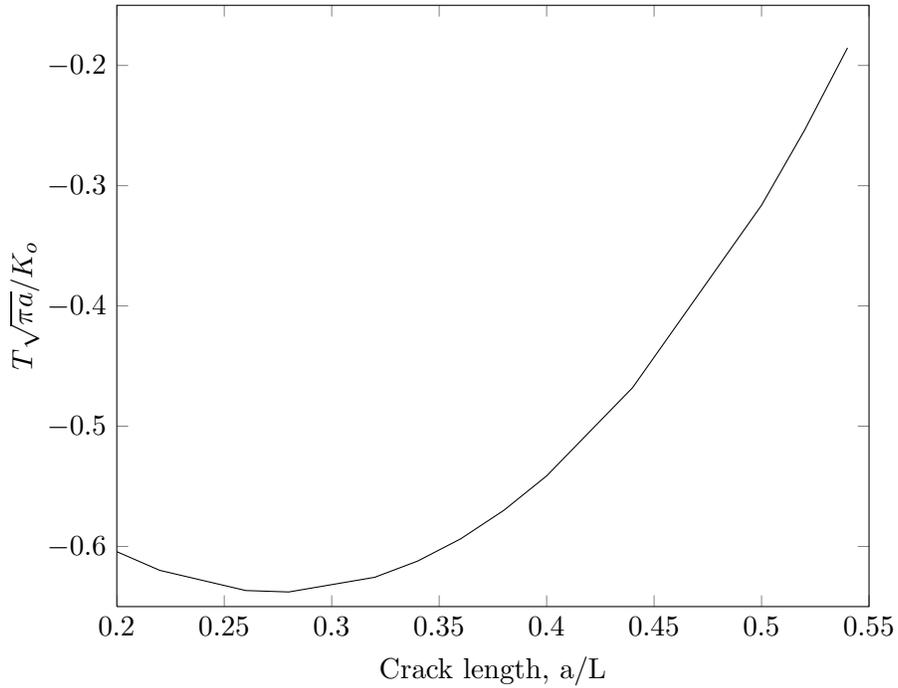
\begin{figure}[htpb]
\centering
\setlength\figureheight{8cm} 
\setlength\figurewidth{10cm}
\subfigure[Stress intensity factor $K_o=\sqrt{K_I^2+K_{II}^2}$]{\scalebox{1}{
%
%
%
%
\begin{tikzpicture}

\begin{axis}[%
width=\figurewidth,
height=\figureheight,
scale only axis,
xmin=0.2,
xmax=0.55,
xlabel={Crack length, a/L},
ymin=1.4,
ymax=3.4,
ylabel={$K_o = \sqrt{K_I^2 + K_{II}^2}$}
]
\addplot [
color=black,
solid,
forget plot
]
table[row sep=crcr]{
0.2 1.40734646008255\\
0.22 1.45020924659444\\
0.26 1.55277529382098\\
0.28 1.61270946373623\\
0.32 1.75129172275365\\
0.34 1.83081204261333\\
0.36 1.91790940059144\\
0.38 2.01329731001185\\
0.4 2.11782272393013\\
0.44 2.35846872896807\\
0.5 2.81961129319421\\
0.52 3.00766928443574\\
0.54 3.21721496643452\\
};
\end{axis}
\end{tikzpicture}%
}}
\subfigure[T-stress $T \sqrt{\pi a}/K_o$]{\scalebox{1}{
%
%
%
%
\begin{tikzpicture}

\begin{axis}[%
width=\figurewidth,
height=\figureheight,
scale only axis,
xmin=0.2,
xmax=0.55,
xlabel={Crack length, a/L},
ymin=-0.65,
ymax=-0.15,
ylabel={$T \sqrt{\pi a}/K_o$}
]
\addplot [
color=black,
solid,
forget plot
]
table[row sep=crcr]{
0.2 -0.604343921867707\\
0.22 -0.61987508886716\\
0.26 -0.636645571372909\\
0.28 -0.637838379848372\\
0.32 -0.625645672090421\\
0.34 -0.612156772243302\\
0.36 -0.593656827700072\\
0.38 -0.570072898407393\\
0.4 -0.541322110141923\\
0.44 -0.467929041743222\\
0.5 -0.316229648322418\\
0.52 -0.25391555513525\\
0.54 -0.185370556726417\\
};
\end{axis}
\end{tikzpicture}%
}}
\caption{Crack growth along the interface: the variation of stress intensity factor and T-stress with the increasing crack length, $a/L$.}
\label{fig:bimatcrkgwth}
\end{figure}

\paragraph{Crack growth along the interface}
In this case, we assume that the crack grows along the interface. Consider a bimaterial plate with an edge crack subjected to uniform tension (see \fref{fig:crkbimat}). The initial crack length $a/L$ is assumed to be 0.2 and the crack increment at each step is assumed to be 0.2. A structured quadrilateral mesh (51$\times$ 102) with 5 layers of elements replaced by the SBFEM is used for this study. Like the XFEM, as the crack advances, a new region is identified which is intersected by the discontinuous surface and instead of enriching with additional enrichment functions, the SBFEM is employed to compute the stiffness matrix and for further post-processing of the results. \fref{fig:bimatcrkmesh} shows the meshes at the initial step and at intermediate steps. The influence of the crack propagation on the numerical stress intensity factors are shown in \fref{fig:bimatcrkgwth}. It is observed that with increasing crack length $a/L$, the stress intensity factor and the T-stress increases.

\subsubsection{Bimaterial plate with a center crack in tension} 
\begin{figure}[htpb]
\centering
\scalebox{0.6}{\input{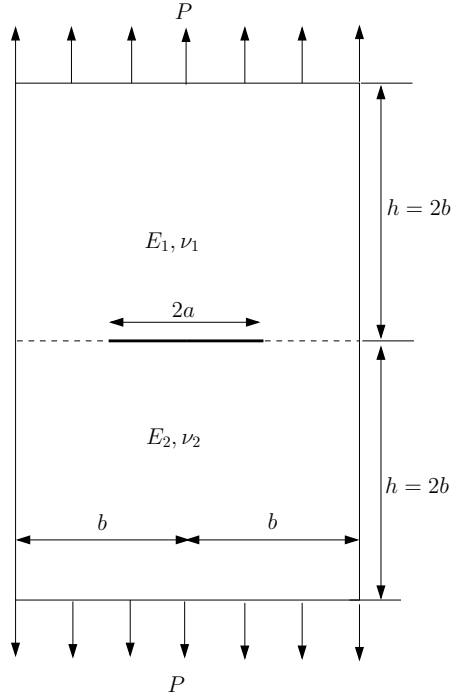}}
\caption{Bi-material plate with a center crack in tension: geometry and boundary conditions.}
\label{fig:crkbimatcenter}
\end{figure}

Consider a bimaterial plate with a center crack subjected to uniform far field tension. The geometry, loading and boundary conditions are shown in \fref{fig:crkbimatcenter}. The material is assumed to be isotropic with Young's modulus $E_1$ and $E_2$ and Poisson's ratio $\nu_1$ and $\nu_2$. Due to symmetry, only one half of the plate is considered. A state of plane strain is considered in this example. Based on the observation from the previous example, a structured quadrilateral mesh (51 $\times$ 102 elements) with 5 layers of elements were replaced with the SBFEM in the vicinity of the crack tip. The stress intensity factors $(K_{\rm I}, K_{\rm II})$ and the T-stress are computed from their definitions (see Section \ref{gsiftress}).

\begin{table}[htpb]
\centering
\renewcommand{\arraystretch}{1.2}
\caption{Bimaterial plate with a center crack in tension: Stress intensity factors and the T-stress, where $K_o = \sqrt{K_I^2 + K_{II}^2}$.}
\begin{tabular}{lrrrrrrr}
\hline 
$E_1/E_2$ & \multicolumn{2}{c}{$K_{\rm I}/\sqrt{\pi a}$} & \multicolumn{2}{c}{$K_{\rm II}/\sqrt{\pi a}$} & \multicolumn{3}{c}{$T\sqrt{\pi a}/K_o$} \\
\cline{2-8}
 & Ref.~\cite{song2005} & Present & Ref.~\cite{song2005} & Present  & Ref.~\cite{song2005} & Ref.~\cite{sladeksladek1997} & Present  \\
\hline
1	& 1.189 &	1.1893 & 	0.000 &	0.0000 &	-1.06 & -1.07 & -1.0552 \\
2	& 1.179 &	1.1798	& -0.055 &	-0.0566 & -0.718 & -0.73 & -0.7144 \\
5	& 1.148 &	1.1483	 & -0.104	 & -0.1053 & -0.379 & -0.39 & -0.3770 \\
10 &	1.123 &	1.1237 & 	-0.123 &	-0.1240 &-0.216	 & -0.23 & -0.2147 \\ 
\hline
\end{tabular}
\label{table:tstressbimaterialinterfacecrk}
\end{table}

\begin{figure}
\centering
\setlength\figureheight{8cm} 
\setlength\figurewidth{10cm}
\subfigure[]{
%
%
%
%
\begin{tikzpicture}

\begin{axis}[%
width=\figurewidth,
height=\figureheight,
scale only axis,
xmin=-200,
xmax=200,
xlabel={$\text{Angle }\theta\text{ (deg)}$},
ymin=-0.5,
ymax=2,
ylabel={Normalized singular stresses},
legend style={at={(0.766964285714285,0.751190476190477)},anchor=south west,draw=black,fill=white,legend cell align=left}
]
\addplot [
color=black,
only marks,
mark=o,
mark options={solid}
]
table[row sep=crcr]{
-174.289406862499 0.13440142037433\\
-163.300755766006 0.32146266822819\\
-153.434948822922 0.470816112311738\\
-145.007979801441 0.571850899411192\\
-138.012787504183 0.632162028787249\\
-131.987212495817 0.668098367023369\\
-124.992020198559 0.685939186364133\\
-116.565051177078 0.6811727956669\\
-106.699244233994 0.649946524594062\\
-95.7105931374997 0.604426196606496\\
-84.2894068625004 0.57511594629339\\
-73.3007557660064 0.590560246731945\\
-63.434948822922 0.654094165080597\\
-55.0079798014413 0.747304596924349\\
-48.0127875041833 0.848751479618651\\
-41.9872124958167 0.951009962537113\\
-34.9920201985587 1.07796966140385\\
-26.565051177078 1.23330862311094\\
-16.6992442339936 1.39848406456519\\
-5.71059313749964 1.5333434369158\\
5.71059313749964 0.434518651717176\\
16.6992442339936 0.428655255092781\\
26.565051177078 0.388817629942456\\
34.9920201985587 0.335480601094625\\
41.9872124958167 0.283539360482602\\
48.0127875041833 0.237768537659602\\
55.0079798014413 0.186978935347778\\
63.434948822922 0.136604993218954\\
73.3007557660064 0.101442433043431\\
84.2894068625003 0.0997741853440532\\
95.7105931374996 0.138932221259797\\
106.699244233994 0.204250393505581\\
116.565051177078 0.269792393642601\\
124.992020198559 0.31777929337442\\
131.987212495817 0.344173420094373\\
138.012787504183 0.352173497310141\\
145.007979801441 0.342658581112779\\
153.434948822922 0.300065563687248\\
163.300755766006 0.206674697362167\\
174.289406862499 0.0554028398856065\\
};
\addlegendentry{$\sigma{}_{\text{xx}}$};

\addplot [
color=black,
solid
]
table[row sep=crcr]{
-174.289406862499 -0.00124332607063266\\
-163.300755766006 0.000723082334091067\\
-153.434948822922 0.0243062165431545\\
-145.007979801441 0.0724668833690262\\
-138.012787504183 0.136469291913085\\
-131.987212495817 0.2126133150992\\
-124.992020198559 0.318673547780566\\
-116.565051177078 0.471339375425909\\
-106.699244233994 0.672703265321389\\
-95.7105931374997 0.902175074899032\\
-84.2894068625004 1.11432533399339\\
-73.3007557660064 1.26467015047593\\
-63.434948822922 1.33992765830793\\
-55.0079798014413 1.3567038996494\\
-48.0127875041833 1.33995907965191\\
-41.9872124958167 1.30837981934769\\
-34.9920201985587 1.25373260140197\\
-26.565051177078 1.17565364856358\\
-16.6992442339936 1.08626725700781\\
-5.71059313749964 1.01584591172382\\
5.71059313749964 0.99620108745791\\
16.6992442339936 1.01172915853516\\
26.565051177078 1.04869747294852\\
34.9920201985587 1.09072194989766\\
41.9872124958167 1.12696247871382\\
48.0127875041833 1.15438351032209\\
55.0079798014413 1.17858119753935\\
63.434948822922 1.18904252675572\\
73.3007557660064 1.16648087638836\\
84.2894068625003 1.08992285230096\\
95.7105931374996 0.954135863046108\\
106.699244233994 0.781551933603464\\
116.565051177078 0.608252719470311\\
124.992020198559 0.460030766680879\\
131.987212495817 0.344950805554519\\
138.012787504183 0.253375017516175\\
145.007979801441 0.165538300565929\\
153.434948822922 0.0851842401639464\\
163.300755766006 0.0274120208365745\\
174.289406862499 0.00224624403213558\\
};
\addlegendentry{$\sigma{}_{\text{yy}}$};

\addplot [
color=black,
dashed
]
table[row sep=crcr]{
-174.289406862499 -0.00581765500550634\\
-163.300755766006 0.0182611313254285\\
-153.434948822922 0.081213449766426\\
-145.007979801441 0.156845641356792\\
-138.012787504183 0.226927977980405\\
-131.987212495817 0.28370451992205\\
-124.992020198559 0.346826394844296\\
-116.565051177078 0.404184100857882\\
-106.699244233994 0.433423798319488\\
-95.7105931374997 0.407972069353667\\
-84.2894068625004 0.317992709012348\\
-73.3007557660064 0.186051685453518\\
-63.434948822922 0.0515198180186466\\
-55.0079798014413 -0.0576737061727596\\
-48.0127875041833 -0.13371641249859\\
-41.9872124958167 -0.1813890391155\\
-34.9920201985587 -0.215451040615939\\
-26.565051177078 -0.218443342521008\\
-16.6992442339936 -0.168474224325441\\
-5.71059313749964 -0.0552464441704165\\
5.71059313749964 0.0416874287482724\\
16.6992442339936 0.0835448294660288\\
26.565051177078 0.107135172791979\\
34.9920201985587 0.107890573110134\\
41.9872124958167 0.091515059865881\\
48.0127875041833 0.0657107706779571\\
55.0079798014413 0.0198475611106207\\
63.434948822922 -0.0534623734464033\\
73.3007557660064 -0.155871766334553\\
84.2894068625003 -0.274057791557087\\
95.7105931374996 -0.37828935251935\\
106.699244233994 -0.440216374324567\\
116.565051177078 -0.453295867133945\\
124.992020198559 -0.430920545159713\\
131.987212495817 -0.391004494342452\\
138.012787504183 -0.346610964828797\\
145.007979801441 -0.281813996003284\\
153.434948822922 -0.197130226434493\\
163.300755766006 -0.103280130494254\\
174.289406862499 -0.0257329135220018\\
};
\addlegendentry{$\tau{}_{\text{xy}}$};

\end{axis}
\end{tikzpicture}%
}
\subfigure[]{\input{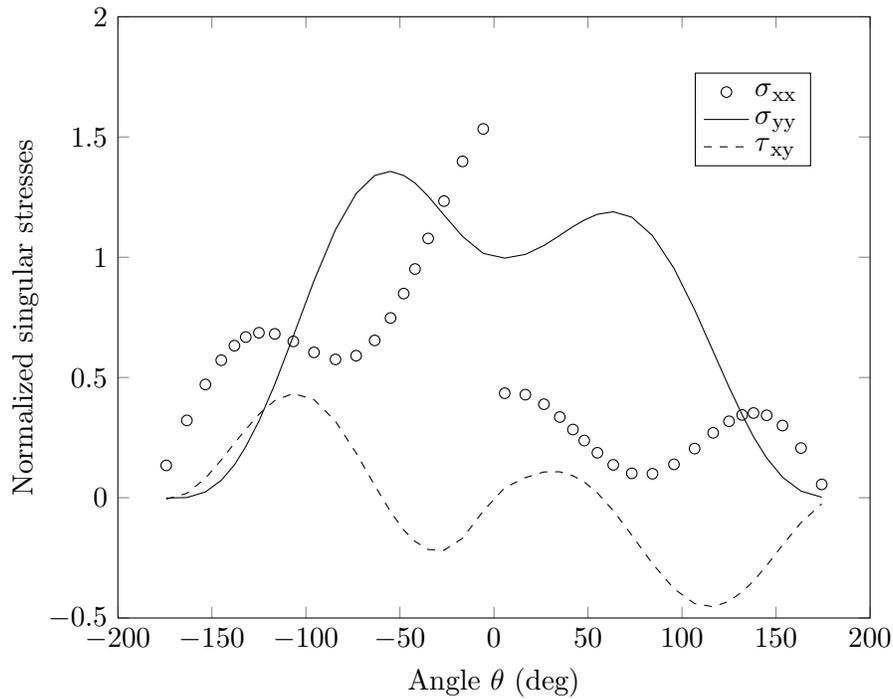}
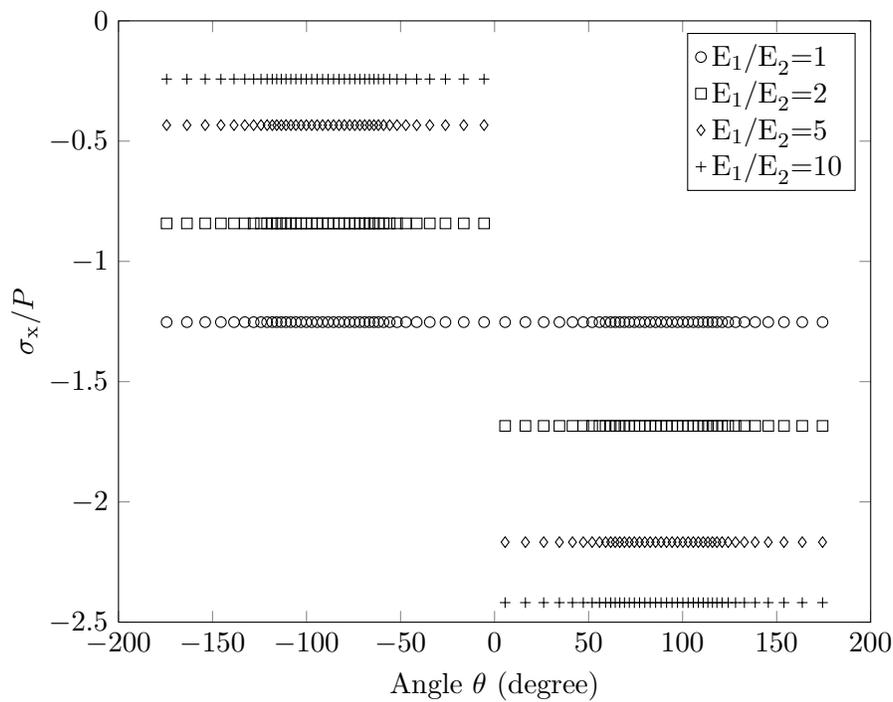}
\caption{Bimaterial plate with a center crack in tension: (a) Singular stress distribution for $E_1/E_2=5$ and (b) the T-stress distribution for various material combinations.}
\label{fig:tstrResultsbimatintecrk}
\end{figure}

Table \ref{table:tstressbimaterialinterfacecrk} presents the stress intensity factors ($K_{\rm I}$ \& $K_{\rm II}$) and the T-stress for various combinations of material property: $E_1/E_2=$ 1, 2, 5 and 10 for a crack length $a/b=$ 0.5. The Poisson's ratio is assumed to be constant with $\nu_1=\nu_2=$ 0.3. When $E_1/E_2 \neq 1$, the order of the singularity is complex and the stress intensity factors are computed from \Eref{eq:GSIFdef-case2-3-1} with the characteristic length $L=2a$. The numerically computed stress intensity factors are normalized with $P\sqrt{\pi a}$. The T-stress is evaluated by using \Eref{eq:TStress-1} to the side of Material 2 (at $\theta= 0^-$) and is normalized with $K_o/\sqrt{\pi a}$. The results from the present formulation are compared with the boundary element technique~\cite{sladeksladek1997} and with the scaled boundary formulation~\cite{song2005} and a very good agreement is observed. \fref{fig:tstrResultsbimatintecrk} shows the singular stress distribution around the crack tip for Young's modulus ratio $E_1/E_2=$ 5 and the T-stress distribution for various combinations of material property. The singular stresses $\bvsig^s$ and the T-stress are evaluated by using \Eref{eq:singularStress_theta0} \&  \Eref{eq:TStress-1}, respectively. It is seen that, as the ratio $E_1/E_2$ increases, the T-stress increases.

\subsubsection{Bimaterial strip with edge loads} 
\begin{figure}[htpb]
\centering
\subfigure[]{\scalebox{0.7}{\input{./Figures/bimatstrip.pstex_t}}}
\subfigure[]{\includegraphics[scale=0.40]{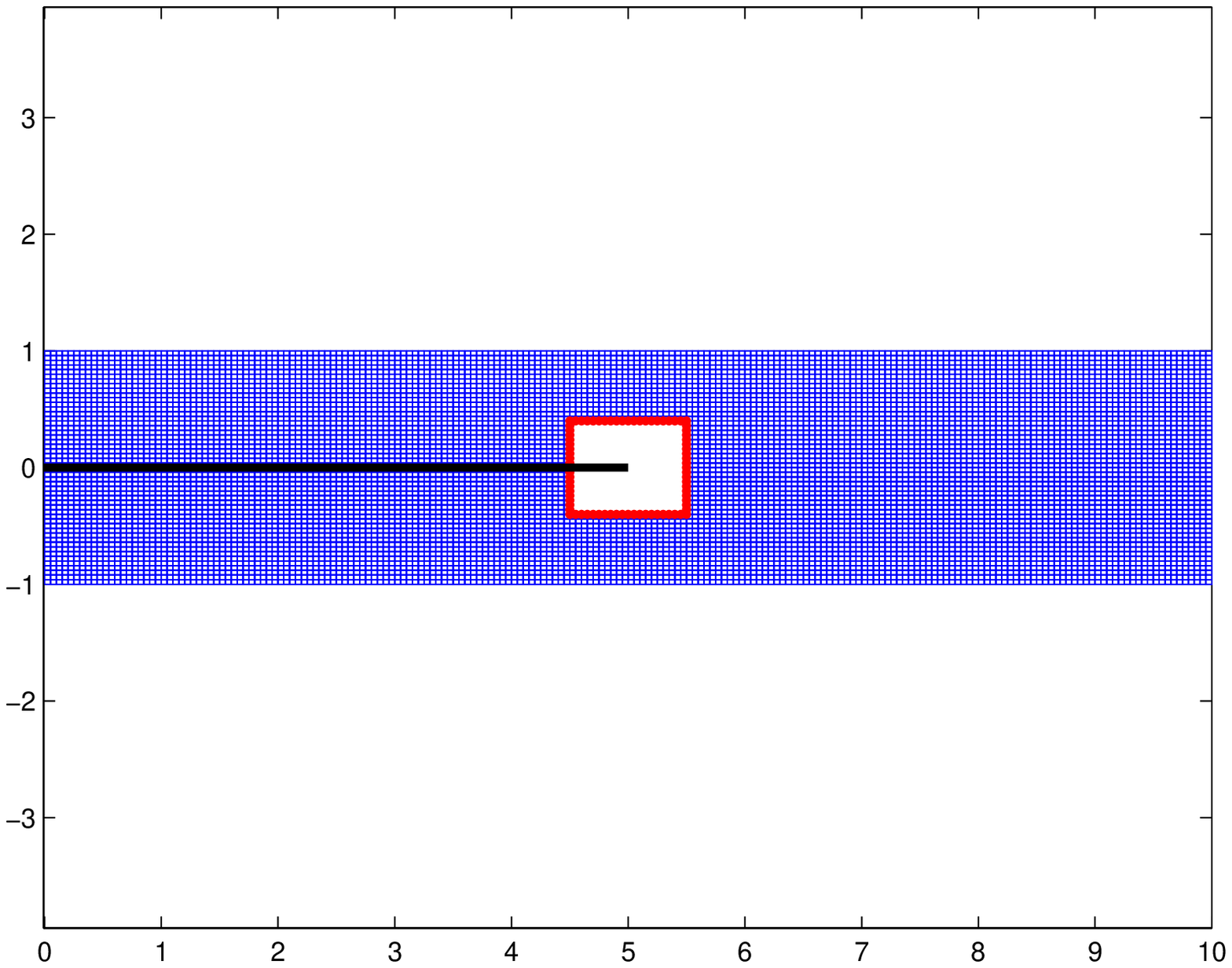}}
\subfigure[]{\includegraphics[scale=0.40]{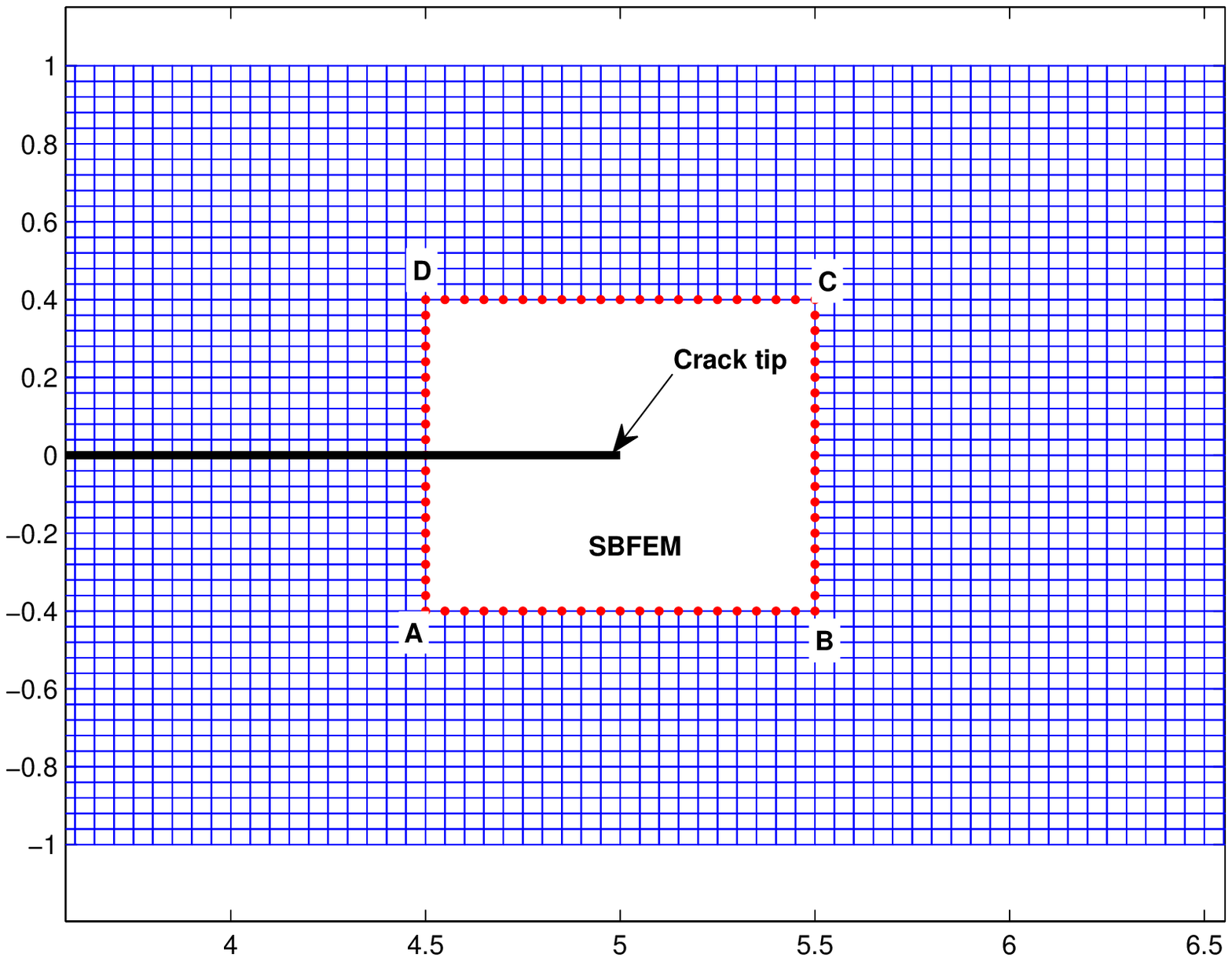}}
\caption{Bi-material strip under edge loads: (a) geometry and boundary conditions; (b) and (c) FE discretization and zoomed in view. The region $A-B-C-D$ is represented by the SBFEM.}
\label{fig:bimstripproblem}
\end{figure}

In this example, consider a bimaterial strip constrained at the right end while a concentrated load $P_o=$ 1 applied at the left end as shown in \fref{fig:bimstripproblem}. The strip contains two layers of thickness $h_1$ and $h_2$ and is assumed to be in a state of plane strain condition. The length of the bi-material strip is $L=10h_2$ and $h_2=$1 is considered for this study. A structured mesh of 201 $\times$ 51 quadrilateral elements are considered with 10 layers of elements replaced with the scaled boundary domain. The analytical expression for the T-stress is given in~\cite{kimvlassak2006}. Table \ref{table:tstressbimaterialstrip} shows the normalized T-stress computed from the current approach for various combinations of material parameters, where $T_1$ and $T_2$ are the T-stress evaluated at the side of the material 1 (at $\theta=0^+$) and the material 2 (at $\theta=0^-$), respectively  . It can be seen that the results from the present formulation are in excellent agreement with the results available in the literature. It is emphasized that the present formulation does not require a priori knowledge of the asymptotic fields.

\begin{table}[htpb]
\centering
\caption{Bimaterial strip with edge loads: material constants and the normalized T-stress. The length of the bimaterial strip, $L=$10$h_2$, the crack length $a=L/2$ and concentrated load, $P_o=$ 1.}
\begin{tabular}{lrrrrrrrrr}
\hline 
Case & $E_1/E_2$ & $\nu_1$ & $\nu_2$ & \multicolumn{2}{c}{Ref.~\cite{kimvlassak2006}} & \multicolumn{2}{c}{Ref.~\cite{yuwu2012}} & \multicolumn{2}{c}{Present} \\
\cline{5-10}
   &  &  &  & T$_1^{\rm ref}$/$\sigma_o$ & T$_2^{\rm ref}$/$\sigma_o$ & T$_1^{\rm ref}$/$\sigma_o$ & T$_2^{\rm ref}$/$\sigma_o$ & T$_1^{\rm ref}$/$\sigma_o$ & T$_2^{\rm ref}$/$\sigma_o$ \\
\hline
1 & 7/3 & 1/3 & 1/3 & 0.0709 & 0.0304 & 0.0702 & 0.0301 & 0.07012 & 0.03005\\
 2& 20/9 & 1/4 & 1/8 & 0.0784 & 0.0336 & 0.0773 & 0.0331 & 0.07698 & 0.03299 \\
 3 & 4 & 2/5 & 2/5 & 0.1310 & 0.0328 & 0.1317 & 0.0329 & 0.13150 & 0.03288 \\
 4 & 4 & 1/4 & 1/4 & 0.1424 & 0.0356 & 0.1410 & 0.0353 & 0.14081 & 0.03520 \\
\hline
\end{tabular}
\label{table:tstressbimaterialstrip}
\end{table}

\subsection{Crack terminating at an interface}
\begin{figure}
\centering
\scalebox{0.6}{\input{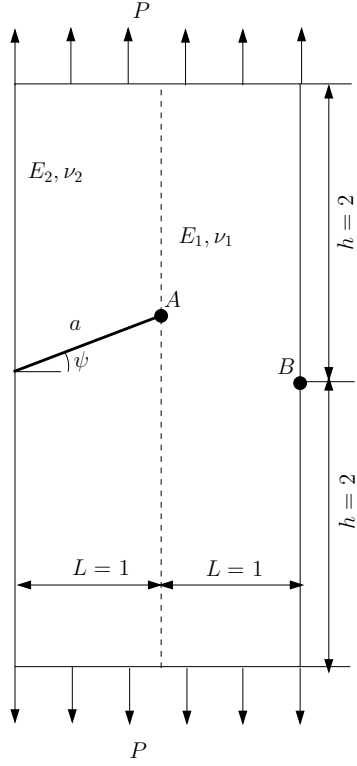}}
\caption{Crack terminating at the material interface: geometry and boundary condition. The crack is assumed to be in material 1 and terminating at the interface between the two materials. $\psi$ is the measured counter clockwise from the positive $x$ axis and determines the inclination of the crack to the interface.}
\label{fig:crackatinter}
\end{figure}

\begin{figure}
\centering
\setlength\figureheight{8cm} 
\setlength\figurewidth{10cm}
\input{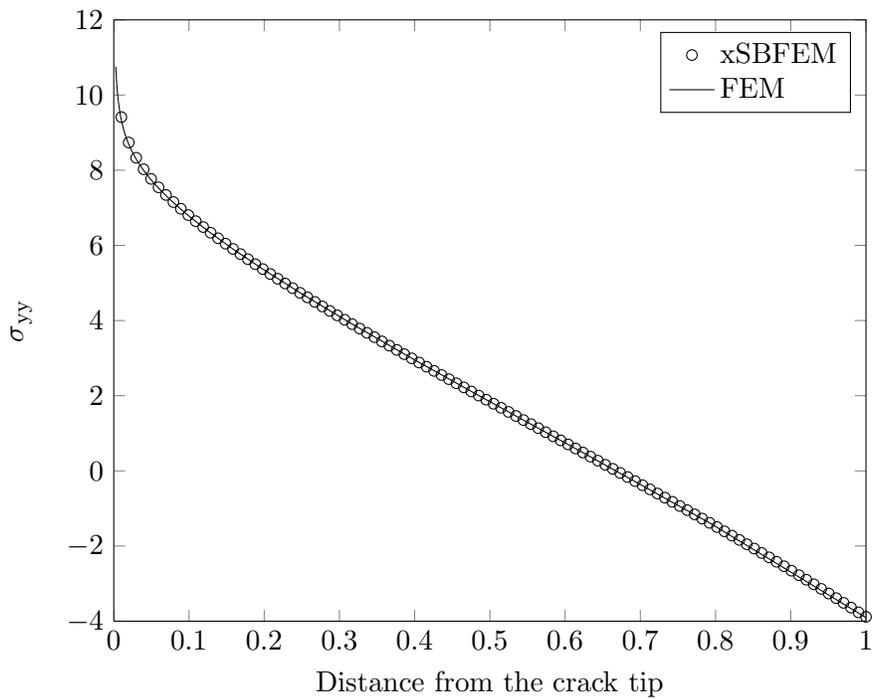}
\caption{Crack terminating at the material interface: stress $\sigma_{yy}$ ahead of the crack tip (along the path $A-B$, see \fref{fig:crackatinter}). The results are compared with the finite element analysis in which a conforming mesh is generated. The ratio of Young's modulus is $E_1/E_2=$ 1000 with $E_2=$ 1.}
\label{fig:tstrResults}
\end{figure}

\begin{figure}
\centering
\setlength\figureheight{8cm} 
\setlength\figurewidth{10cm}
\subfigure[$E_1/E_2=$ 100]{\input{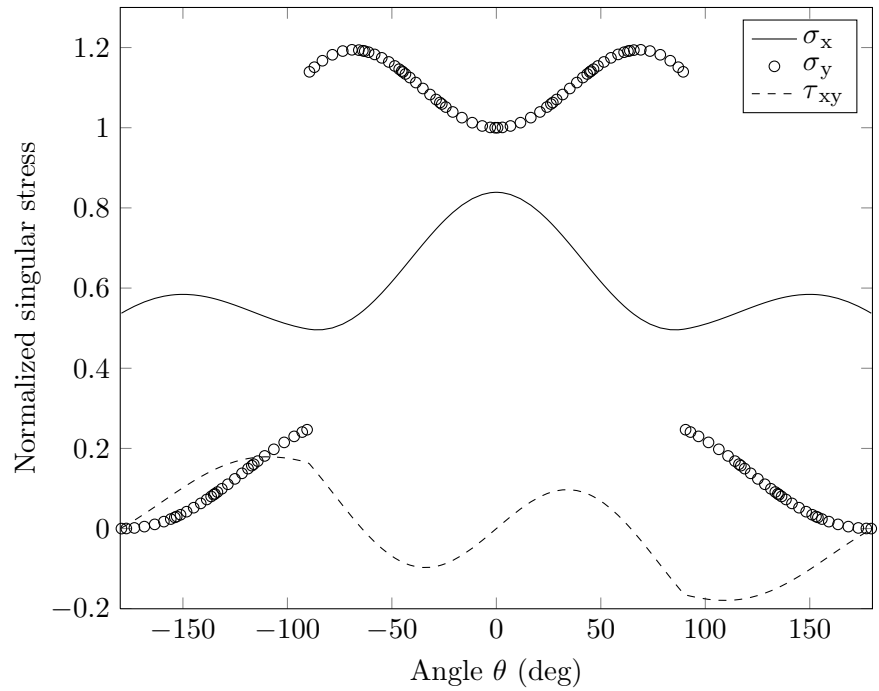}}
\subfigure[]{\input{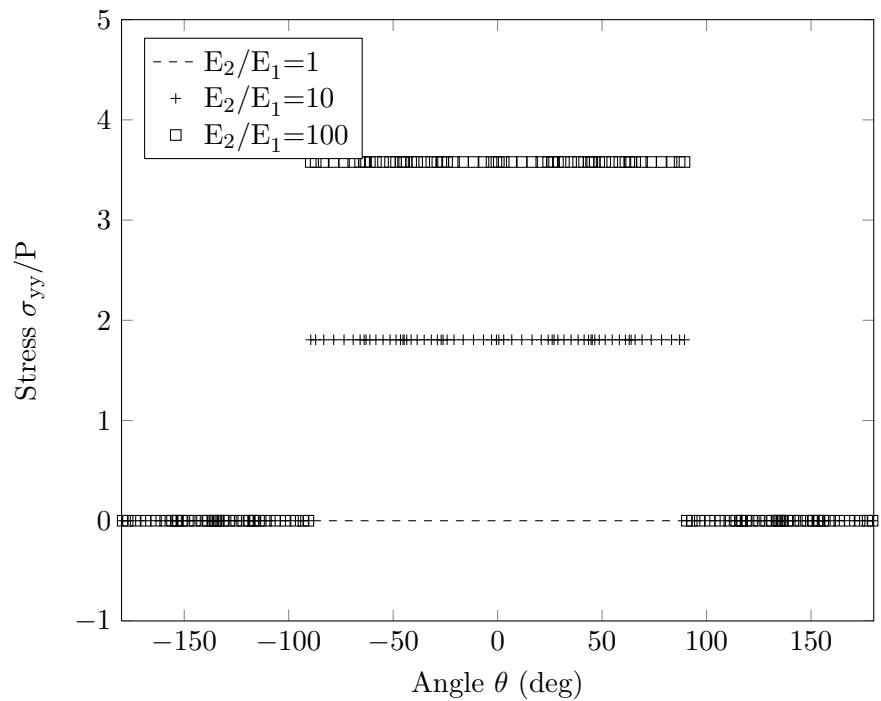}}
\caption{Angular distribution of singular stresses and T-stress for various combination of material properties for a crack terminating at the material interface.}
\label{fig:angdistriTcrkatInt}
\end{figure}

\begin{figure}
\centering
\setlength\figureheight{8cm} 
\setlength\figurewidth{10cm}
\input{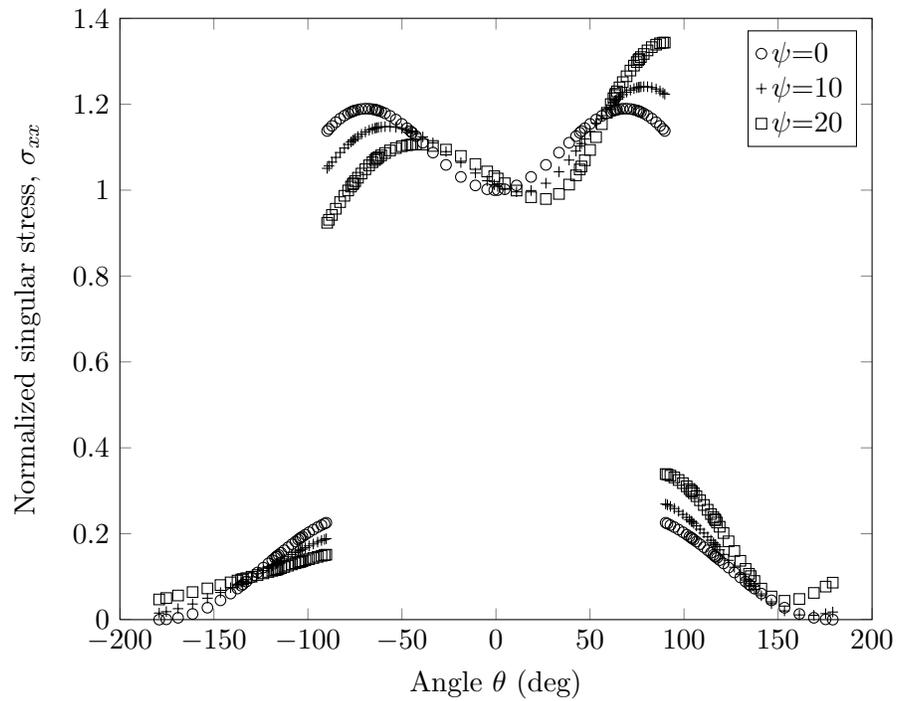}
\caption{Angular distribution of singular stresses and T-stress for various combination of material properties for a crack terminating at the material interface for various crack orientation. The ratio of Young's modulus is $E_1/E_2=$ 100 with $E_2=$ 1.}
\label{fig:angdistriTcrkatIntvaricrkangle}
\end{figure}

In this last example, consider a crack terminating at the bimaterial interface. The crack is assumed to be in material 2 and terminating at the interface between the two materials (see \fref{fig:crackatinter}) and the plate is subjected to a far field tension $P$. Appropriate Dirichlet boundary conditions are enforced to arrest any rigid body motion. A structured quadrilateral mesh $(50 \times 100)$ with 5 layers of element replaced with the SBFEM is used for this study. \fref{fig:tstrResults} shows the stress in the $y-$direction $(\sigma_{yy}$ ahead of the crack tip). The results from the present formulation is compared with the stresses computed using a finite element model with a conforming mesh. \fref{fig:angdistriTcrkatInt} shows the angular stress distribution ahead of the crack tip for $E_1/E_2=$ 100 and the T-stress distribution for various combination of material properties when the crack is perpendicular to the interface, i.e., $\psi=$0$^\circ$. It is seen that, as the ratio $E_2/E_1$ increases, the T-stress increases. The angular distribution of $\sigma_{yy}$ for various crack orientation $\psi$ is shown in \fref{fig:angdistriTcrkatIntvaricrkangle}. When the crack is perpendicular to the interface, the stress distribution is symmetric, whilst in any other case it is unsymmetric. It is observed that the crack orientation strongly influences the stress distribution ahead of the crack tip. 
 
\begin{figure}
\centering
\setlength\figureheight{8cm} 
\setlength\figurewidth{10cm}
%
%
%
%
\begin{tikzpicture}

\begin{axis}[%
width=\figurewidth,
height=\figureheight,
scale only axis,
xmin=0.4,
xmax=0.8,
xlabel={a/L},
ymin=-1,
ymax=8,
ylabel={Mode I stress intensity factor, $K_I$},
legend style={draw=black,fill=white,legend cell align=left},
legend pos = north east
]
\addplot [
color=black,
solid
]
table[row sep=crcr]{
0.43 2.52284\\
0.44922552168068 2.80774\\
0.468434692380206 3.0678\\
0.487662907011386 3.30305\\
0.506893634225566 3.59337\\
0.526124350054968 3.90043\\
0.545355077695916 4.2409\\
0.564585809561905 4.62105\\
0.583816558883783 5.0478\\
0.603047316932814 5.52722\\
};
\addlegendentry{$\text{K}_\text{I}\text{ (}\psi\text{=20)}$};

\addplot [
color=black,
solid,
mark=square,
mark options={solid}
]
table[row sep=crcr]{
0.43 2.60383\\
0.449230769194463 2.83488\\
0.468461538403685 3.07294\\
0.487692307618425 3.3038\\
0.506923076841262 3.59254\\
0.526153846066876 3.89853\\
0.545384615294521 4.23882\\
0.564615384523588 4.61891\\
0.583846153754149 5.04556\\
0.603076922984914 5.52705\\
0.62230769221527 6.07362\\
0.641538461444379 6.69229\\
0.660769230674566 7.41024\\
};
\addlegendentry{$\text{K}_\text{I}\text{ (}\psi\text{=0)}$};

\end{axis}

\begin{axis}[%
width=\figurewidth,
height=\figureheight,
scale only axis,
xmin=0.4,
xmax=0.8,
xminorticks=true,
ymin=-0.15,
ymax=0.5,
xmode=log,
hide x axis,
axis y line*=right,
ylabel={Mode II stress intensity factor, $K_{II}$},
ylabel near ticks,
legend style={draw=black,fill=white,legend cell align=left},
legend pos = south east
]

\addplot [
color=black,
dashed
]
table[row sep=crcr]{
0.43 0.39831\\
0.44922552168068 -0.09963\\
0.468434692380206 0.04773\\
0.487662907011386 0.02347\\
0.506893634225566 0.00799\\
0.526124350054968 -0.00054\\
0.545355077695916 -0.00023\\
0.564585809561905 -0.00123\\
0.583816558883783 -0.00091\\
0.603047316932814 -0.00014\\
};
\addlegendentry{$\text{K}_{\text{II}}\text{ (}\psi\text{=20)}$};

\addplot [
color=black,
dashed,
mark=o,
mark options={solid}
]
table[row sep=crcr]{
0.43 -8e-05\\
0.449230769194463 2e-05\\
0.468461538403685 1e-05\\
0.487692307618425 2e-05\\
0.506923076841262 1e-05\\
0.526153846066876 1e-05\\
0.545384615294521 1e-05\\
0.564615384523588 2e-05\\
0.583846153754149 1e-05\\
0.603076922984914 2e-05\\
0.62230769221527 2e-05\\
0.641538461444379 -7e-05\\
0.660769230674566 2e-05\\
};
\addlegendentry{$\text{K}_{\text{II}}\text{ (}\psi\text{=0)}$};

\end{axis}
\end{tikzpicture}%
\caption{Mode I and II stress intensity factor variation as a function of crack advancement. The crack increment is set in advance and the direction of the crack propagation is computed from \Eref{eqn:probdirection}. }
\label{fig:sifvariationcrkgwth}
\end{figure}
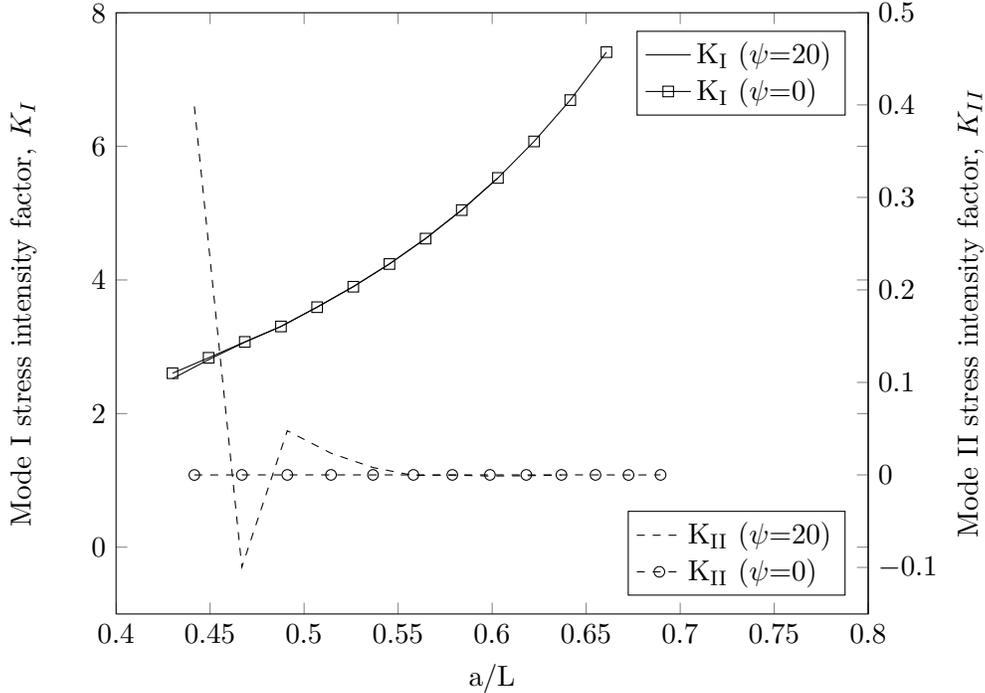

Next we consider the case of a crack deflecting into the material. A crack impinging the interface between two dissimilar materials may advance by either propagating along the interface or deflecting into the material. The competition between these two scenarios have been thoroughly dealt in the literature~\cite{hehutchinson1989,hehutchinson1989a,hutchinsonsuo1992,zhangsuo2007}. The material properties are assumed to be $E_2=100E_1$. We assume that the initial crack is in material 2 and has deflected into the material 1. In this example, material 1 is homogeneous and isotropic and hence the crack growth is governed by the maximum hoop stress criterion, which states that the crack will propagate from its tip in the direction $\theta_c$ where the circumferential stress $\sigma_{\theta \theta}$ is maximum. The critical angle is computed by solving the following equation:
\begin{equation}
K_I \sin \theta_c + K_{II} (3 \cos \theta_c - 1) = 0
\label{eqn:cirstress}
\end{equation}
Solving \Eref{eqn:cirstress} gives the crack propagation angle:
\begin{equation}
\theta_c = 2 \arctan \left[ \frac{-2\left( \frac{K_{II}}{K_I} \right)}{ 1+\sqrt{1 + 8\left( \frac{K_{II}}{K_I} \right)^2}}  \right]
\label{eqn:probdirection}
\end{equation}
The plate is subjected to a far field tension $P$ (see \fref{fig:crackatinter}). A structured quadrilateral mesh (51$\times$102) with 5 layers of element replaced with the SBFEM is employed in this example. \fref{fig:sifvariationcrkgwth} shows the variation of the stress intensity factors (mode I and mode II) with crack length $a/L$. The influence of the initial orientation of the crack is also shown. Since the loading is mode I, we expect the crack to grow in a straight line. For a crack impinging at an angle $\psi$, initially $K_{II}$ is not zero (see \fref{fig:sifvariationcrkgwth}) and hence the crack advances at an angle that is not equal to zero. However, due to the far field tension, as the crack advances, $K_{II}$ becomes zero and the crack propagates in a straight line. If conventional XFEM was employed, a separate set of enrichment functions should be used when the crack tip lies on the interface and these functions are to be replaced with asymptotic fields for homogeneous material when the crack deflects into the material interface. It is noted that in the extended scaled boundary method, no such modification is necessary and moreover the stress intensity factors and the T-stress can be computed from the definitions. 

\begin{figure}[htpb]
\centering
\includegraphics[scale=0.6]{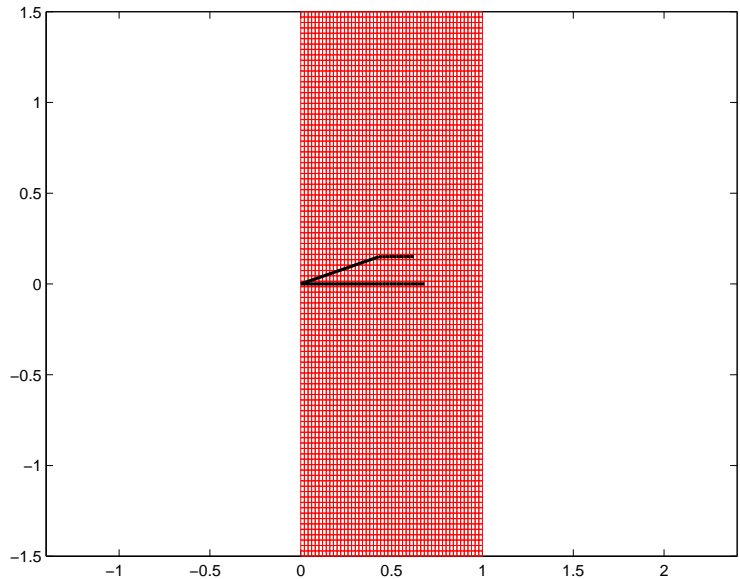}
\caption{Crack trajectory for a crack impinging at the interface and deflected into material 1. Two cases are considered: (a) a crack perpendicular to the interface and (b) crack oriented at an angle. The solid black line denotes the discontinuous surface.}
\label{fig:cracktrajectory}
\end{figure}

\section{Conclusion}
In this paper, we studied the strength of the singularity of cracks aligned to bimaterial interfaces and cracks terminating at the interface by employing a semi-analytical approach. Later, the stress intensity factors and the T-stress was computed by using a combined XFEM and a semi-analytical approach. The salient features of this technique are: (i) no a priori knowledge of the asymptotic fields is required and (ii) the stiffness matrix of the region enriched for a particular crack is obtained directly without a special numerical integration technique. This eliminates the need to compute the enrichment functions as in~\cite{duartekim2008,menkbordas2010,mousavigrinspun2011,zhu2012}, whilst employing the Heaviside functions allows the discontinuous surfaces to be presented independent of the underlying mesh.The crack growth along the interface and the deflection of the crack into a material is also discussed. The proposed technique can readily be coupled with the existing XFEM/GFEM code and can be applied to study fracture behaviour in heterogeneous materials.

\section*{Acknowledgements} 
S Natarajan would like to acknowledge the financial support of the School of Civil and Environmental Engineering, The University of New South Wales for his research fellowship for the period September 2012 onwards. 


\section*{References}
\bibliographystyle{elsarticle-num}
\bibliography{SinguStress}

\begin{thebibliography}{10}
\expandafter\ifx\csname url\endcsname\relax
  \def\url#1{\texttt{#1}}\fi
\expandafter\ifx\csname urlprefix\endcsname\relax\def\urlprefix{URL }\fi
\expandafter\ifx\csname href\endcsname\relax
  \def\href#1#2{#2} \def\path#1{#1}\fi

\bibitem{natarajansong2013}
S.~Natarajan, C.~Song, Representation of signular fields without asymptotic
  enrichment in the extended finite element method, International Journal for
  Numerical Methods in Engineering 96 (2013) 813--841.

\bibitem{barsoum1977}
R.~Barsoum, Triangular quater-point elements as elastic and perfectly plastic
  crack tip elements, International Journal for Numerical Methods in
  Engineering 11 (1977) 85--98.

\bibitem{strangfix1973}
G.~Strang, G.~Fix, An analysis of the finite element method, Prentice-Hall,
  1973.

\bibitem{tracey1971}
D.~Tracey, Finite elements for determination of crack-tip elastic stress
  intensity factors, Tech. Rep. Technical Report AD0732837, Army materials and
  mechanics research center, Watertown Mass (1971).

\bibitem{atlurikobayashi1975}
S.~Atluri, A.~Kobayashi, M.~Nakagaki, An assumed displacement hybrid finite
  element model for linear fracture mechanics, International Journal of
  Fracture 11 (1975) 257--271.

\bibitem{benzley1974}
S.~Benzley, Reapplications of singularities with isoparametric finite elements,
  International Journal for Numerical Methods in Engineering 8 (1974) 537--545.

\bibitem{rabczukgracie2010}
T.~Rabczuk, R.~Gracie, J.~Song, T.~Belytschko, Immersed particle method for
  fluid-structure interaction, International Journal for Numerical Methods in
  Engineering 81 (2010) 48--71.

\bibitem{rabczukbelytschko2006}
T.~Rabczuk, T.~Belytschko, Application of particle methods to static fracture
  of reinforced concrete structures, International Journal of Fracture 137
  (2006) 19--49.

\bibitem{rabczukzi2008a}
T.~Rabczuk, G.~Zi, S.~Bordas, H.~Nguyen-Xuan, A geometrically nonlinear three
  dimensional cohesive crack method for reinforced concrete structures,
  Engineering Fracture Mechanics 75 (2008) 4740--4758.

\bibitem{zirabczuk2007}
G.~Zi, T.~Rabczuk, W.~Wall, Extended meshfree methods without branch enrichment
  for cohesive cracks, Computational Mechanics 40 (2007) 367--382.

\bibitem{bacnguyen-xuan2013}
N.~V. Bac, H.~Nguyen-Xuan, L.~Chen, S.~Bordas, X.~Zhuang, G.~Liu, T.~Rabczuk, A
  phantom-node method with edge based strain smoothing for linear elastic
  fracture mechanics, Journal of Applied Mathematics.

\bibitem{chau-dinhzi2012}
T.~Chau-Dinh, G.~Zi, P.~Lee, J.~Song, T.~Rabczuk, A phantom-node method for
  shellmodels with arbitrary cracks, Computers \& Structures 92--93 (2012)
  242--256.

\bibitem{rabczukbelytschko2007}
T.~Rabczuk, T.~Belytschko, A three-dimensional large deformation meshfree
  method for arbitrary evolving cracks, Computer Methods in Applied Mechanics
  and Engineering 196 (2007) 2777--2799.

\bibitem{zhuangaugarde2012}
X.~Zhuang, C.~Augarde, K.~Mathisen, Fracture modeling using meshless methods
  and level sets in 3{D}: framework and modeling, International Journal for
  Numerical Methods in Engineering 92 (2012) 969--998.

\bibitem{rabczukzi2010}
T.~Rabczuk, G.~Zi, S.~Bordas, H.~Nguyen-Xuan, A simple and robust three
  dimensional cracking-particle method without enrichment, Computer Methods in
  Applied Mechanics and Engineering 199 (2010) 2437--2455.

\bibitem{bordasgzi2008}
S.~Bordas, T.~R. nd~G~Zi, Three-dimensional crack initiation, propagation,
  branching and junction in nonlinear materials by extrinsic discontinuous
  enrichment of meshfree methods without asymptotic enrichment, Engineering
  Fracture Mechanics 75 (2008) 943--960.

\bibitem{rabczukbordas2007}
T.~Rabczuk, S.~Bordas, G.~Zi, A three-dimensional meshfree method for
  continuous multiple crack initiation, nucleation and propagation in statics
  and dynamics, Computational Mechanics 40 (2007) 473--495.

\bibitem{melenkbabuvska1996}
J.~Melenk, I.~Babu\v{s}ka, The partition of unity finite element method: basic
  theory and applications, Computer Methods in Applied Mechanics and
  Engineering 39 (1996) 289--314.

\bibitem{belytschkoblack1999}
T.~Belytschko, T.~Black, Elastic crack growth in finite elements with minimal
  remeshing, International Journal for Numerical Methods in Engineering 45
  (1999) 601--620.

\bibitem{garciafancello2000}
O.~Garcia, E.~Fancello, C.~Barcellos, C.~Duarte, {$hp$-clouds in Mindlin's
  thick plate model}, International Journal for Numerical Methods in
  Engineering 47 (2000) 1381--1400.

\bibitem{rabczukbelytschko2004}
T.~Rabczuk, T.~Belytschko, Cracking particles: a simplified meshfree method for
  arbitrary evolving cracks, International Journal for Numerical Methods in
  Engineering 61 (2004) 2316--2343.

\bibitem{rabczukzi2008}
T.~Rabczuk, G.~Zi, A.~Gerstenberger, W.~Wall, A new crack tip element for the
  phantom-node method with arbitrary cohesive cracks, International Journal for
  Numerical Methods in Engineering 75 (2008) 577--599.

\bibitem{nguyen-xuanliu2013}
H.~Nguyen-Xuan, G.~Liu, S.~Bordas, S.~Natarajan, T.~Rabczuk, An adaptive
  singular {ES-FEM} for mechanics problems with singular fields of arbitrary
  order, Computer Methods in Applied Mechanics and Engineering 253 (2013)
  252--273.

\bibitem{belytschkoorgan1995}
T.~Belytschko, D.~Organ, Y.~Krongauz, A coupled finite element - element-free
  {Galerkin} method, Computational Mechanics 17 (1995) 186--195.

\bibitem{bird2012}
G.~E. Bird, The coupled dual boundary element-scaled boundary finite element
  method for efficient fracture mechanics, Ph.D. thesis, Durham University
  (2012).

\bibitem{duartekim2008}
C.~Duarte, D.-J. Kim, Analysis and applications of a generalized finite element
  method with global-local enrichment functions, Computer Methods in Applied
  Mechanics and Engineering 197 (2008) 487--504.

\bibitem{waismanbelytschko2008}
H.~Waisman, T.~Belytschko, Parametric enrichment adaptivity by the extended
  finite element method, International Journal for Numerical Methods in
  Engineering 73 (2008) 1671--1692.

\bibitem{menkbordas2010}
A.~Menk, S.~Bordas, Numerically determined enrichment functions for the
  extended finite element method and applications to bi-material anisotropic
  fracture and polycrystals, International Journal for Numerical Methods in
  Engineering 83 (2010) 805--828.

\bibitem{mousavigrinspun2011}
S.~Mousavi, E.~Grinspun, N.~Sukumar, Harmonic enrichment functions: A unified
  treatment of multiple intersecting and branched cracks in the extended finite
  element method, International Journal for Numerical Methods in Engineering 85
  (2011) 1306--1322.

\bibitem{mousavigrinspun2011a}
S.~Mousavi, E.~Grinspun, N.~Sukumar, Higher-order extended finite elements with
  harmonic enrichment functions for complex crack problems, International
  Journal for Numerical Methods in Engineering 86 (2011) 560--574.

\bibitem{zhu2012}
Q.-Z. Zhu, On enrichment functions in the extended finite element method,
  International Journal for Numerical Methods in Engineering 91 (2012)
  186--217.

\bibitem{xiaokarihaloo2007}
Q.~Xiao, B.~Karihaloo, Implementationa of hybrid crack element on a general
  finite element mesh and in combination with {XFEM}, Computer Methods in
  Applied Mechanics and Engineering 196 (2007) 1864--1873.

\bibitem{chahinelaborde2008}
E.~Chahine, P.~Laborde, Y.~Renard, {Spider XFEM, an extended finite element
  variant for partially unknown crack-tip displacement}, European Journal of
  Computational Mechanics 17 (2008) 625--636.

\bibitem{rethoreroux2010}
J.~R\'ethore, S.~Roux, F.~Hild, Hybrid analytical and extended finite element
  method {(HAX-FEM)}: {A new enrichment procedure for cracked solids},
  International Journal for Numerical Methods in Engineering 81 (2010)
  269--285.

\bibitem{babuvskabanerjee2012}
I.~Babu\v{s}ka, U.~Banerjee, Stable generalized finite element method
  {(SSGFEM)}, Computer Methods in Applied Mechanics and Engineering 201--204
  (2012) 91--11.

\bibitem{asadpouremohammadi2007}
A.~Asadpoure, S.~Mohammadi, A new approach to simulate the crack with the
  extended finite element method in orthotropic media, International Journal
  for Numerical Methods in Engineering 69 (2007) 2150--2172.

\bibitem{asharimohammadi2011}
S.~E. Ashari, S.~Mohammadi, Delamination analysis of composites by new
  orthotropic bimaterial extended finite element method, International Journal
  for Numerical Methods in Engineering 86 (2011) 1507--1543.

\bibitem{hattorirojas-diaz2012}
G.~Hattori, R.~Rojas-D\'iaz, A.~S\'aez, N.~Sukumar, F.~Garcia-S\'anchz, New
  anisotropic crack-tip enrichment functions for the extended finite element
  method, Computational Mechanics 50 (2012) 591--601.

\bibitem{dolbowgosz2002}
J.~Dolbow, M.~Gosz, On the computation of mixed-mode stress intensity factors
  in functionally graded materials, International Journal of Solids and
  Structures 39 (2002) 2557--2574.

\bibitem{sukumarhuang2004}
N.~Sukumar, Z.~Huang, J.-H. Pr\'evost, Z.~Suo, Partition of unity enrichment
  for bimaterial interface cracks, International Journal for Numerical Methods
  in Engineering 59 (2004) 1075--1102.

\bibitem{liuxiao2004}
X.~Liu, Q.~Xiao, B.~Karihaloo, {XFEM for direct evaluation of mixed mode SIFs
  in homogeneous and bi-material}, International Journal for Numerical Methods
  in Engineering 59 (2004) 1103--1118.

\bibitem{bouhalashao2013}
L.~Bouhala, Q.~Shao, Y.~Koutsawa, A.~Younes, P.~N. {n}ez, A.~Makradi,
  S.~Belouettar, An {XFEM} crack-tip enrichment for a crack terminating at a
  bi-material interface, Engineering Fracture Mechanics 102 (2013) 51--64.

\bibitem{kimpaulino2003}
J.~Kim, G.~Paulino, The interaction integral for fracture of orthotropic
  functionally graded materials: evaluation of stress intensity factors,
  International Journal of Solids and Structures 40 (2003) 3967--4001.

\bibitem{courtingardin2005}
S.~Courtin, C.~Gardin, C.~B\'ezine, H.~Hamouda, {Advantages of the J-integral
  approach for calculating stress intensity factors when using the commercial
  finite element software ABAQUS}, Engineering Fracture Mechanics 72 (2005)
  2174--2185.

\bibitem{bergerkarageorghis2007}
J.~Berger, A.~Karageorghis, P.~Martin, Stress intensity factor computation
  using the method of fundamental solutions: Mixed-mode problems, International
  Journal for Numerical Methods in Engineering 69 (2007) 469--483.

\bibitem{morais2007}
A.~de~Morais, Calproblems of stress intensity factors by the force method,
  Engineering Fracture Mechanics 74 (2007) 739--750.

\bibitem{passieuxgravouil2011}
J.~Passieux, A.~Gravouil, J.~R\'ethor\'e, {Direct estimation of generalized
  stress intensity factors using a three-scale concurrent multigrid X-FEM},
  International Journal for Numerical Methods in Engineering 85 (2011)
  1648--1666.

\bibitem{wolfsong2001}
J.~Wolf, C.~Song, The scaled boundary finite-element method - a fundamental
  solution-less boundary element method, Computer Methods in Applied Mechanics
  and Engineering 190 (2001) 5551--5568.

\bibitem{deekswolf2002}
A.~Deeks, J.~Wolf, A virtual work derivation of the scaled boundary finite
  element method for elastostatics, Computational Mechanics 28 (2002) 489--594.

\bibitem{ooisong2012}
E.~T. Ooi, C.~Song, F.~Tin-Loi, Z.~Yang, Polygon scaled boundary finite
  elements for crack propagation modeling, International Journal for Numerical
  Methods in Engineering 91 (2012) 319--342.

\bibitem{natarajanooi2013}
S.~Natarajan, E.~T. Ooi, I.~Chiong, C.~Song, {Convergence and accuracy of
  ddisplacement based finite element formulations over arbitrary polygons:
  Laplace interpolants, strain smoothing and scaled boundary polygon
  formulation}, Finite Elements in Analysis and Design~Accepted in press.

\bibitem{williams1957}
M.~Williams, On the stress distribution at the base of a stationary crack,
  {Journal of Applied Mechanics, ASME} 24 (1957) 109--114.

\bibitem{song2005}
C.~Song, {Evaluation of power-logarithmic singularities, T-stress and higher
  order terms of in-plane singular stress field at cracks and multi-material
  corners}, Engineering Fracture Mechanics 72 (2005) 1498--1530.

\bibitem{yuwu2012}
H.~Yu, L.~Wu, H.~Li, {T-stress eevaluation of an interface crack in the
  materials with complex interfaces}, International Journal of Fracture 177
  (2012) 25--37.

\bibitem{sladeksladek1997}
J.~Sladek, V.~Sladek, {Evaluations of the T-stress for interface cracks by the
  boundary element method}, Engineering Fracture Mechanics 56 (1997) 813--825.

\bibitem{kimvlassak2006}
J.~Kim, J.~Vlassak, {T-stress of a bi-material strip under generalized edge
  loads}, International Journal of Fracture 142 (2006) 315--322.

\bibitem{matsumtotanaka2000}
T.~Matsumto, M.~Tanaka, R.~Obara, Computation of stress intensity factors of
  interface cracks based on interaction energy release rates and {BEM}
  sensitivity analysis, Engineering Fracture Mechanics 65 (2000) 683--702.

\bibitem{hehutchinson1989}
M.-Y. He, J.~W. Hutchinson, Crack deflection at an interface between dissimilar
  elastic materials, International Journal of Solids and Structures 25 (1989)
  1053--1067.

\bibitem{hehutchinson1989a}
M.-Y. He, J.~W. Hutchinson, materials a crack out of an interface, Journal of
  Applied Mechanics 56 (1989) 270--278.

\bibitem{hutchinsonsuo1992}
J.~W. Hutchinson, Z.~Suo, Mixed mode cracking in layered materials, Adv. Appl.
  Mech. 29 (1992) 63--191.

\bibitem{zhangsuo2007}
Z.~Zhang, Z.~Suo, Split singularities and the competition between crack
  penetration and debond at a bimaterial interface, International Journal of
  Solids and Structures 44 (2007) 4559--4573.

\end{thebibliography}

\end{document}